# A New Continuum-Based Thick Shell Finite Element for Soft Biological Tissues in Dynamics: Part 1 - Preliminary Benchmarking Using Classic Verification Experiments

Bahareh Momenan, Ph.D.
*Department of Mechanical Engineering, University of Ottawa,*
*161 Louis Pasteur, Ottawa, K1N 6N5 Ontario, Canada*
*baharehmomenan@gmail.com*

Michel R. Labrosse, Ph.D.
*Department of Mechanical Engineering, University of Ottawa,*
*161 Louis Pasteur,* Colonel By Hall A-213, *Ottawa, K1N 6N5 Ontario, Canada*
*labrosse@eng.uottawa.ca*

For the finite element simulation of thin soft biological tissues in dynamics, shell elements, compared to volume elements, can capture the whole tissue thickness at once, and feature larger critical time steps. However, the capabilities of existing shell elements to account for irregular geometries, and hyperelastic, anisotropic 3D deformations characteristic of soft tissues are still limited.

As improvement, we developed a new general nonlinear thick continuum-based (CB) shell finite element (FE) based on the Mindlin-Reissner shell theory, with large bending, large distortion and large strain capabilities, embedded in the updated Lagrangian formulation and explicit time integration. We performed numerical benchmark experiments available from the literature that focus on engineering linear elastic materials, which, verified and proved the new thick CB shell FE to: 1) be accurate an efficient 2) be powerful in handling large 3D deformations, curved geometries, 3) accommodate coarse distorted meshes, and 4) achieve comparatively fast computational times. The new element was also insensitive to three types of locking (shear, membrane and volumetric), and warping effects. The capabilities of the present thick CB shell FE in the biomedical realm are illustrated in a companion article (Part 2), in which anisotropic incompressible hyperelastic constitutive relations are implemented and verified.

## 1. Introduction

Understanding and simulation of the mechanical behaviour of living tissues and organs are two important, interrelated goals towards which finite element (FE) methods have been favoured. For bulky organs, such as the liver and the brain, volume finite elements like tetrahedrals and bricks available in current commercial finite element (FE) programs[1–3] are effective; however, they are less well suited for the analysis of blood vessels and heart valves in that several elements through the thickness of these organs are required to properly account for bending and/or through thickness effects. In dynamic applications, where the critical time steps used for time integration are proportional to the size of the smallest volume elements, this translates into comparatively long computational times. On the other hand, shell elements have the potential to capture the whole thickness of such structures at once, and are known to feature larger critical time steps than volume elements. Therefore, they represent a promising avenue for the FE simulation of relatively thin organs under complex loadings. Although many shell elements already exist in the literature in general, and some shell elements have been proposed for the analysis of soft biological tissues (e.g [4,5]), their capabilities may be limited in modelling irregular and complex geometries, and hyperelastic, large, anisotropic 3D deformations characteristic of soft tissues. The present work aims to address this issue.

Continuum-based (CB) shell finite elements, as developed herein, are derived directly from continuum elements by imposing the kinematic and/or kinetic assumptions of the shell theory of interest [6], with independent translational and rotational degrees of freedom [7]. When based on the Kirchhoff-Love theory for the linear analysis of thin shells, the normal strain (due to large in-plane stretching) and the transverse shear strains (due to rotation of the shell fibers) are not considered [8–12]. Note that in shell theory language, a fiber is the normal to the reference surface of the shell in the undeformed configuration, and is totally unrelated to the collagen fibers commonly present in soft biological tissues. Although large in-plane stretching is formulated in [13] and [14], the normal strains are not calculated, and nor is the thickness of the shell element updated, because of the fiber inextensibility assumption. Thus, volume evolution according to the constitutive relations is ignored. Such limitations arise from the use of only one coordinate system (CS) whose normal is corotational with the fibers, as dictated by the Kirchhoff-Love theory, and in which the zero normal stress is enforced [8–14]. In the mixed interpolation of tensorial components (MITC), the thickness of the element measured along the fiber remains constant during the deformation (hence only small strains are considered), and the accuracy with which transverse shear stresses are predicted depends to a great degree on the mesh used and the geometric distortions of the element [15–17].

Although improved three-dimensional MITC shell elements have been proposed, featuring five (or six) degrees of freedom at each node plus two (or three) additional degrees of freedom to represent the thickness straining and wrapping of the fibers [18–20], another approach to handle large deformations has been based on the implementation of the modified Mindlin-Reissner theory. In such implementation, two independent coordinate systems (CSs) have been introduced to handle the kinematic and the kinetic constraints [21–23], and to allow for large in-plane stretching and rotations of the fibers without creating artificial stiffening, thereby preventing shear locking. However, the thickness update has only been proposed in the context of hypoelasticity (rate form constitutive relations) [21], and application of this formulation to hyperelastic materials such as biological soft tissues would require extensive time integrations that would only add complexity.

Therefore, in developing the present thick CB shell FE, we aimed to combine the advantages and circumvent the limitations of the existing CB shell elements. Specifically, we employed two independent CSs (namely lamina and fiber) to implement the kinematic and kinetic assumptions of the modified Mindlin-Reissner shell theory. Independency of the lamina and fiber CSs permits large rotations of the fibers and inclusion of transverse shears. To allow for complex geometries and large 3D deformations, we enabled large rotations of the independent lamina and fiber CSs at every node in the element. This was achieved by expressing the measures of deformation (Jacobian, deformation gradient, Almansi and Green-Lagrange strains), strain displacement transformation matrices, and constitutive relations, within the lamina CS at each node, as opposed to the global ones in current formulations. We also updated the element thickness such that volume evolution dictated by the linear elastic or hyperelastic constitutive relations was satisfied, as opposed to using rate form or MITC formulations. To minimize the computational cost, while maintaining accuracy with large deformations and rotations, we employed the updated Lagrangian (UL) formulation and an explicit time integration [1–3,6–9,13,15–17,24–30]. Finally, an efficient eigen value approach was implemented to update the critical time step for shell elements with complex geometries and nonlinear material properties.

To verify the accuracy and efficiency characteristics of the new CB shell FE, its large bending and distortion capabilities, and its shear and membrane locking insensitivity, we performed numerical benchmark experiments with different geometries, material properties, and loading conditions. Because most of the detailed verification tests available in the literature involve structures made of engineering materials, the focus herein, as a first step, will be linear and nonlinear isotropic elastic materials. However, in a companion article (Part 2), we developed the kinetics description and the implementation of anisotropic incompressible hyperelastic constitutive relations that enable the CB shell FE to model rubber-like materials and soft biological tissues.

## 2. Methods

### 2.1. *Geometry and kinematics*

We adopted a 9-noded quadrilateral shell element (three nodes per side, and one central node) to represent curved shell geometries and achieve a high degrees of polynomial completeness [31]. The curvilinear parent coordinates $(r, s, t)$ of a 9-noded CB shell element are shown in Fig. 1, left. In this study, the middle surface ($t = 0$, also called lamina) is taken as the reference surface. The continuum representation of the geometry (Fig. 1, right) in any configuration is defined as:

$$^{\beta}y(r,s,t) = \sum_{a=1}^{n_{en}} N_a(r,s)\,{}^{\beta}\bar{y}_a(r,s) + \frac{t}{2}\sum_{a=1}^{n_{en}} N_a(r,s)\;{}^{\beta}h_a\;{}^{\beta}\hat{Y}_a \tag{1}$$

where,

the initial and the current configurations are represented by $\beta = 0$ and $\beta = \tau$ respectively,

$\bar{\bigcirc}$: the accent bar represents the reference surface quantities,

${}^{\beta}\bar{y}_a(r,s)$: is the position vector of a nodal point $a$ in configuration $\beta$,

$N_a(r,s)$: is the 2D shape function associated with node $a$,

$n_{en}$: is the number of element nodes,

${}^{\beta}\hat{Y}_a$: is a unit vector emanating from node $a$ in the fiber direction (also called director).

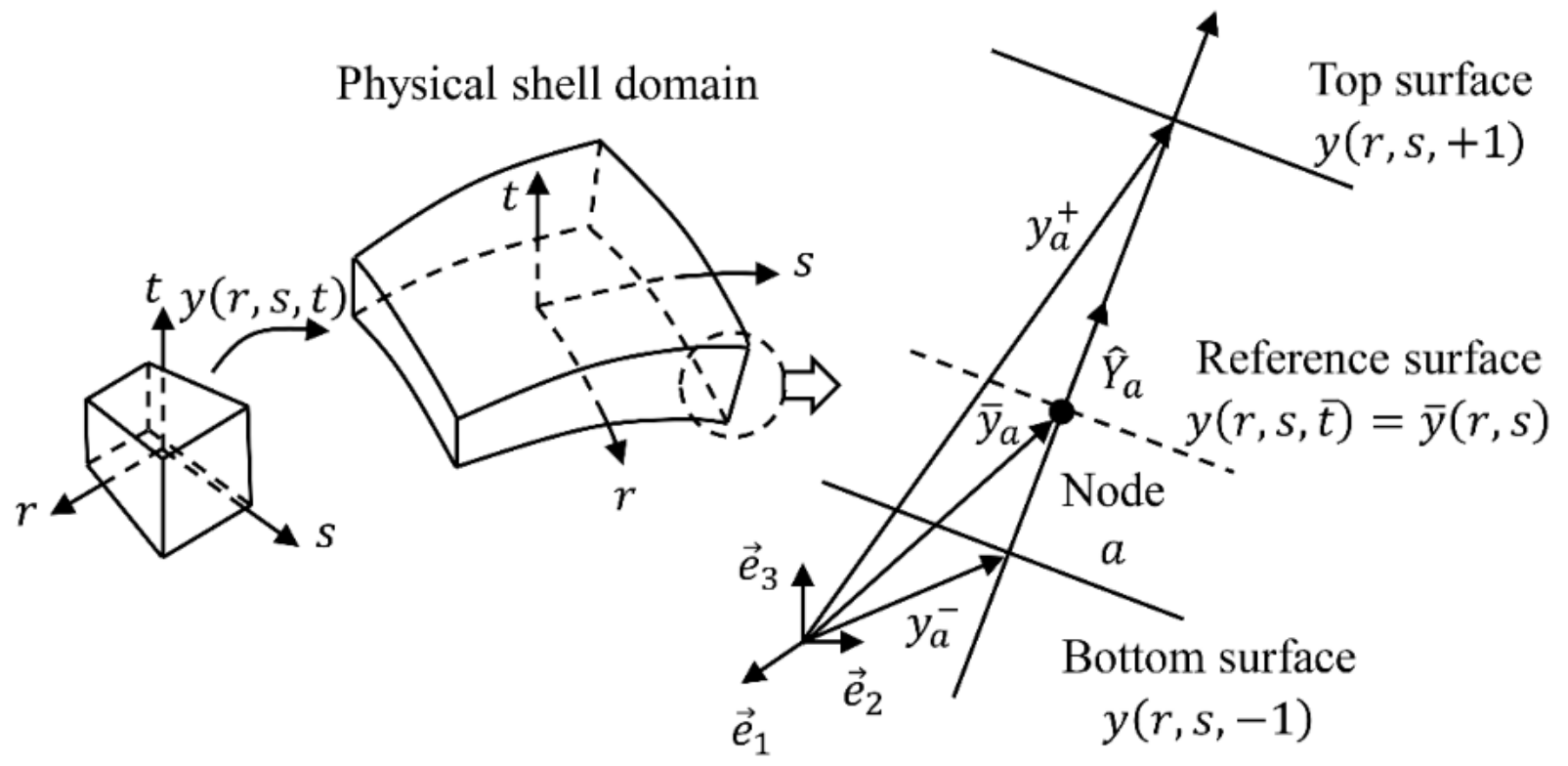

Fig. 1. Left: mapping of a general cube into the physical shell element. Right: illustration of (1) when $\boldsymbol{\beta} = \mathbf{0}$. Please refer to text for details.

In addition, $h_a$ is a nodal fiber length parameter which, in the undeformed configuration ($\beta = 0$), is defined as ${}^0h_a = \|y_a^+ - y_a^-\|$, where $y_a^+$ and $y_a^-$ are the position vectors of the top and bottom surfaces, respectively, of the shell along each nodal fiber (Fig. 1, right).

The kinematic expressions were obtained by replacing the displacement variables with the coordinate variables of (1). This gives:

$$ {}^{\tau+\Delta\tau}_{\beta}u(r,s,t) = \sum_{a=1}^{n_{en}} N_a(r,s)\, {}^{\tau+\Delta\tau}_{\beta}\bar{u}_a(r,s) + \frac{t}{2}\sum_{a=1}^{n_{en}} N_a(r,s)\, {}^{\tau+\Delta\tau}h_a\, {}^{\tau+\Delta\tau}_{\beta}\hat{U}_a \qquad (2) $$

where,

${}^{\tau+\Delta\tau}_{\beta}\bar{u}_a(r,s)$: is the displacement of a node in the reference surface from the reference configuration ($\beta$) to the configuration at time $\tau + \Delta\tau$,

${}^{\tau+\Delta\tau}_{\beta}\hat{U}_a$: is the displacement of a director (in terms of rotations) at a specific node,

${}^{\tau+\Delta\tau}h_a$: is the nodal fiber length in the configuration at time $\tau + \Delta\tau$.

The derivation of the displacement of a director at a specific node ${}^{\tau+\Delta\tau}_{\beta}\hat{U}_a$ is presented in Appendix A. The kinematics described in (2) involve, for each node, three translational degrees of freedom, and two rotational degrees of freedom with respect to the fiber.

Ultimately, the position vector of a generic point of the shell at time $\tau + \Delta\tau$ is obtained by adding the position vector (1) to the displacement vector (2), as ${}^{\tau+\Delta\tau}y(r,s,t) = {}^{\tau}y(r,s,t) + {}^{\tau+\Delta\tau}_{\beta}u(r,s,t)$. This can be written as:

$$ {}^{\tau+\Delta\tau}y(r,s,t) = \sum_{a=1}^{n_{en}} N_a(r,s)\, {}^{\tau+\Delta\tau}\bar{y}_a(r,s) + \frac{t}{2}\sum_{a=1}^{n_{en}} N_a(r,s)\, {}^{\tau+\Delta\tau}h_a\, {}^{\tau+\Delta\tau}\hat{Y}_a \ . \qquad (3) $$

A quick comparison between (1) and (3) suggests that ${}^{\tau+\Delta\tau}h_a\, {}^{\tau+\Delta\tau}\hat{Y}_a = {}^{\tau+\Delta\tau}h_a\, {}^{\tau+\Delta\tau}_{\beta}\hat{U}_a + {}^{\beta}h_a\, {}^{\beta}\hat{Y}_a$ . However, according to Fig. 2, this relation is off by ${}^{\tau+\Delta\tau}_{\beta}h_a\, {}^{\beta}\hat{Y}_a$ , in which ${}^{\tau+\Delta\tau}_{\beta}h_a$ refers to the change in fiber length from the reference configuration $\beta$ to the configuration at time $\tau + \Delta\tau$. This effect is understated in small membrane strain applications, and thus has been neglected in the literature. Instead, herein, we calculate the nodal point fiber lengths at time $\tau + \Delta\tau$, and reevaluate (1) using (${}^{\beta}h_a = {}^{\beta}h_a + {}^{\tau+\Delta\tau}_{\beta}h_a$). The update algorithm will be presented in Section 2.6.

## 2.2. *Coordinate systems and transformation matrices*

As mentioned earlier, in keeping with the modified Mindlin-Reissner theory, we employed two independent (lamina and fiber) CSs to impose the kinetic and the kinematic constraints [21,22]. The lamina CS (Fig. 3, left) is naturally defined by the geometry of the shell, such that two of its bases ($\vec{e}_1^{\,l}$ and $\vec{e}_2^{\,l}$) are tangent to the lamina and the third base ($\vec{e}_3^{\,l}$) is normal to the lamina at each stress storage point (derived in Appendix B). The constitutive relations and the zero normal stress condition are applied within the lamina CS, and $\vec{e}_3^{\,l}$ is not generally tangent to the fiber direction (Fig. 3, right). In contrast, the fiber CS is constructed (see Appendix C) such that one of its bases (denoted as $\vec{e}_3^{\,f}$) coincides and rotates rigidly with the fiber direction ($\hat{Y}$ ) (Fig. 4) at each node. The fiber CS is used as a reference frame for the rotation increments, and to enforce fiber inextensibility.

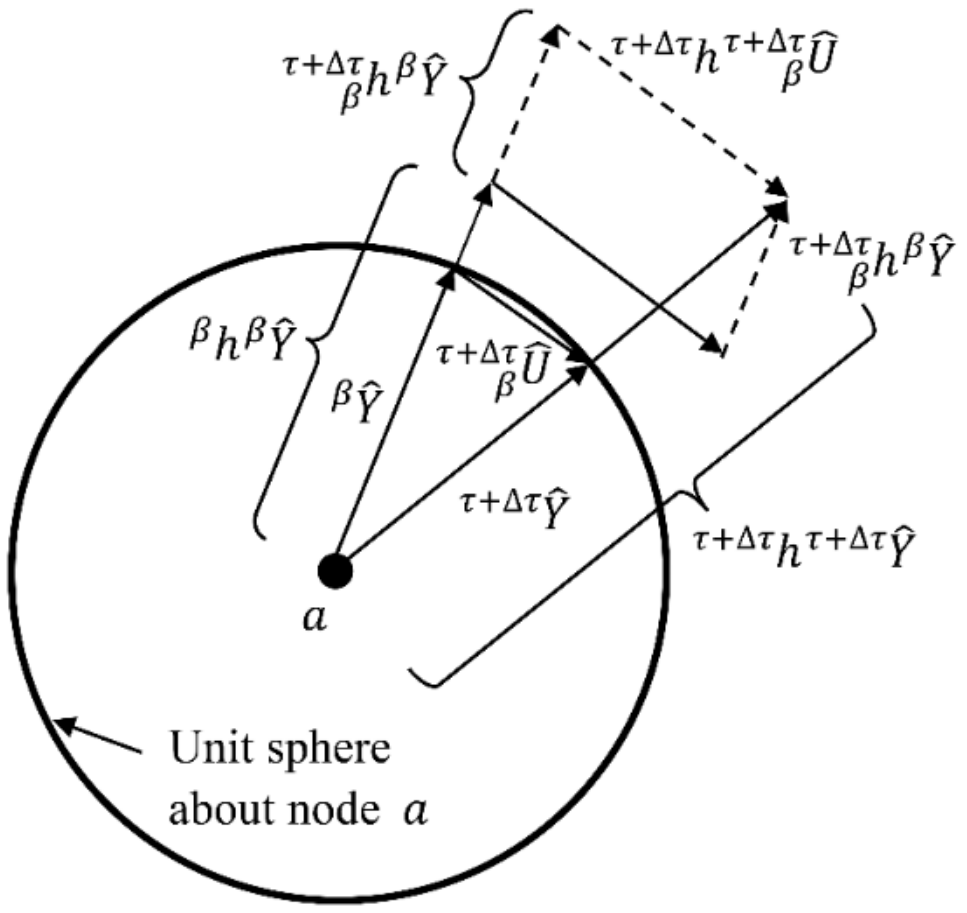

Fig. 2. Kinematics. Right subscripts $\boldsymbol{a}$ denoting the node number are dropped for convenience. Please refer to text for details.

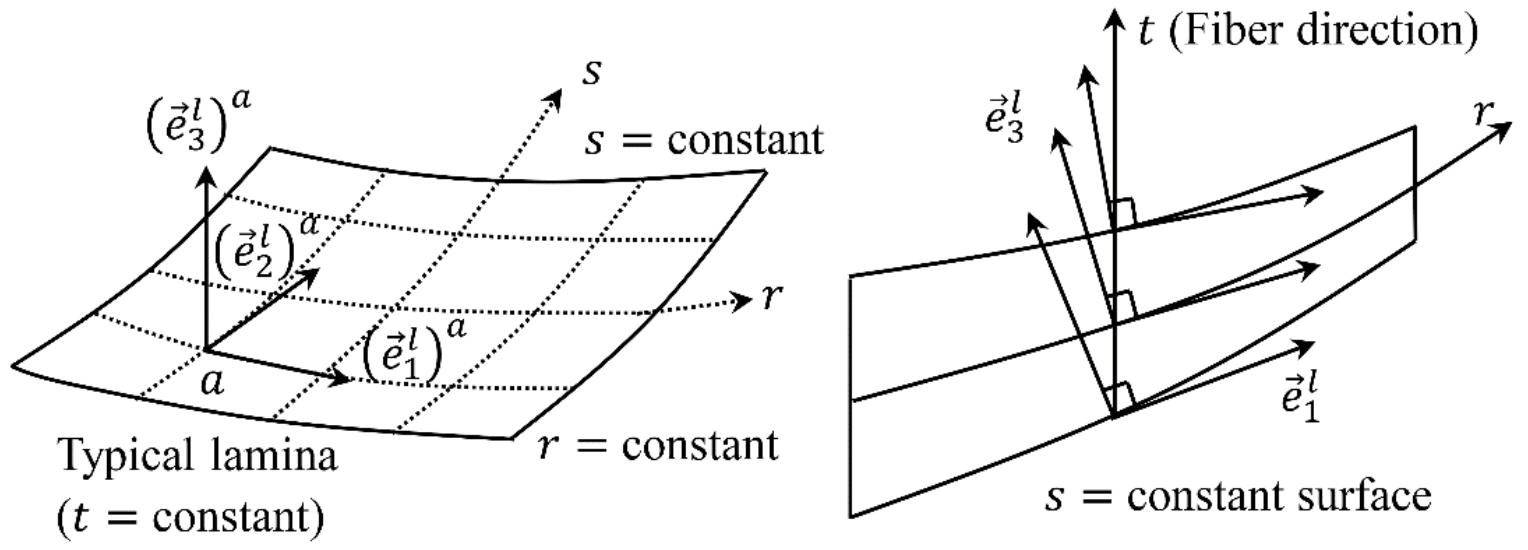

Fig. 3. Left: Lamina CS, at node $\boldsymbol{a}$, shown on a typical lamina. Right: Independence of $\vec{\boldsymbol{e}}_3^{\,l}$ from the fiber direction.

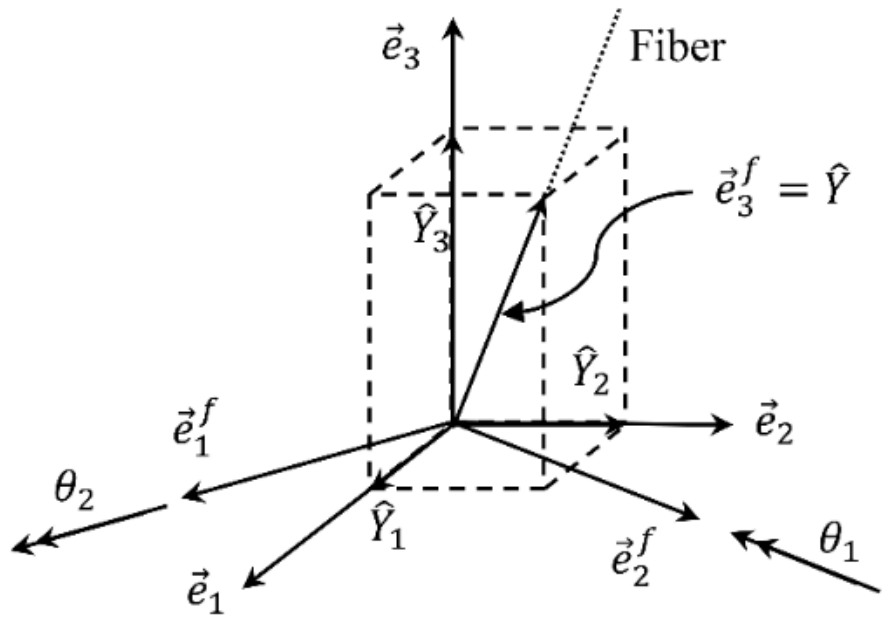

Fig. 4. Illustration of the nodal fiber CS, when, $\vec{\boldsymbol{e}}_3^{\,f}$ is coincident with the fiber direction ($\hat{\boldsymbol{Y}}$ ).

The transformations between the global, fiber and lamina CSs are presented in Fig. 5, where $[s] = [\vec{e}_1^{\,f} \quad \vec{e}_2^{\,f} \quad \vec{e}_3^{\,f}]$, $[q] = [\vec{e}_1^{\,l} \quad \vec{e}_2^{\,l} \quad \vec{e}_3^{\,l}]^T$, and $[r] = [q][s]$. Note that the transformation matrices need to be updated as deformation proceeds. In addition, a subscript can be used to specify the node number (e.g. $[q_a]$ ).

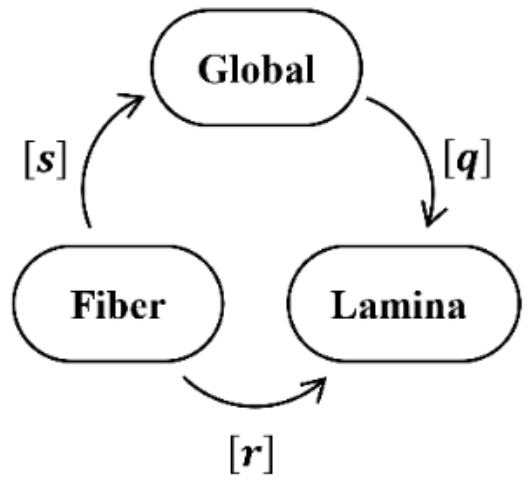

**Fig. 5** Transformation matrices between the three CSs employed in the present work

### 2.3. *Deformation gradient in the lamina coordinate system*

To allow for large discrepancies between the three above CSs at each node, and between nodes within one element, we express the measures of deformation in the lamina. Thus, we obtain the deformation gradient of the current configuration with respect to the reference configuration in the lamina CS from,

$$\left[{}_{\beta}^{\tau}F^{l}\right] = \left[\,{}^{\tau}J^{l}\right]^{T}\left[{}^{\beta}J^{l}\right]^{-T} \tag{4}$$

where, $\left[{}^{\beta}J^{l}\right]$ and $\left[\,{}^{\tau}J^{l}\right]$ are the reference and the current Jacobians, both defined for the nodal lamina CS. Note that the Jacobian at any configuration is the partial derivatives of the position vector in that configuration with respect to the parent CS (for example, see (5) for the definition of the current Jacobian). To derive the Jacobian in the nodal lamina CS, the position vector must be transformed to the nodal lamina CS [22,32] by $\{y^{l}\} = [q]\{y\}$. Based on this principle, we transform the reference and the current Jacobians to the nodal lamina CS by $\left[{}^{\beta}J^{l}\right] = \left[{}^{\beta}J\right]\left[{}^{\beta}q\right]^{T}$ and $\left[\,{}^{\tau}J^{l}\right] = \left[\,{}^{\tau}J\right]\left[\,{}^{\tau}q\right]^{T}$.

$$\left[\,{}^{\tau}J\right] = \begin{bmatrix} \dfrac{\partial\, {}^{\tau}y_1}{\partial r} & \dfrac{\partial\, {}^{\tau}y_2}{\partial r} & \dfrac{\partial\, {}^{\tau}y_3}{\partial r} \\ \dfrac{\partial\, {}^{\tau}y_1}{\partial s} & \dfrac{\partial\, {}^{\tau}y_2}{\partial s} & \dfrac{\partial\, {}^{\tau}y_3}{\partial s} \\ \dfrac{\partial\, {}^{\tau}y_1}{\partial t} & \dfrac{\partial\, {}^{\tau}y_2}{\partial t} & \dfrac{\partial\, {}^{\tau}y_3}{\partial t} \end{bmatrix}. \tag{5}$$

### 2.4. *The UL strain-displacement matrix in the lamina coordinate system*

To derive the kinetics in the current configuration, the displacements equation in the global CS (2) is written in matrix form as $\left\{{}_{\beta}^{\tau}u\right\} = \left[\,{}^{\tau}N\right]\left\{{}_{\beta}^{\tau}u^{a}\right\}$. Expanding, and dropping the left super- and sub-scripts for clarity, one gets:

$$\begin{Bmatrix} u_1 \\ u_2 \\ u_3 \end{Bmatrix} = \begin{bmatrix} \cdots & N_a & 0 & 0 & -\dfrac{(t-\bar{t})}{2}h_a N_a\left(\vec{e}_1^{\,f}\right)_1^a & -\dfrac{(t-\bar{t})}{2}h_a N_a\left(\vec{e}_2^{\,f}\right)_1^a & \cdots \\ \cdots & 0 & N_a & 0 & -\dfrac{(t-\bar{t})}{2}h_a N_a\left(\vec{e}_1^{\,f}\right)_2^a & -\dfrac{(t-\bar{t})}{2}h_a N_a\left(\vec{e}_2^{\,f}\right)_2^a & \cdots \\ \cdots & 0 & 0 & N_a & -\dfrac{(t-\bar{t})}{2}h_a N_a\left(\vec{e}_1^{\,f}\right)_3^a & -\dfrac{(t-\bar{t})}{2}h_a N_a\left(\vec{e}_2^{\,f}\right)_3^a & \cdots \end{bmatrix} \begin{Bmatrix} \vdots \\ u_1^a \\ u_2^a \\ u_3^a \\ \theta_1^a \\ \theta_2^a \\ \vdots \end{Bmatrix}. \tag{6}$$

where $h_a$ refers to the nodal fiber length in the current configuration.

In the UL formulation, the linear strain-displacement transformation matrix $\left[{}_{\tau}^{\tau}B_L\right]$ is obtained by taking the partial derivatives of (6) with respect to the global Cartesian coordinates of the position vector in the current configuration [9,13,25]. However, to have the constitutive relations in the lamina CS, and to allow for large deformations and large strains, we derive the strain-displacement transformation matrix in the lamina CS (denoted by $\left[{}_{\tau}^{\tau}B_L^{l}\right]$) by taking the partial derivatives of the displacement vector in the lamina CS with respect to the current position vector in the lamina CS, as a modification to the formulation of [9,13,25]. As a starting point to the determination of $\left[{}_{\tau}^{\tau}B_L^{l}\right]$, we transform the displacements to the lamina CS $\left\{{}_{\beta}^{\tau}u^{l}\right\}$ by pre-multiplying (6) with the current transformation matrix $\left[\,{}^{\tau}q_a\right]$ evaluated for each node. Then, we take partial derivatives of each row with respect to each of the parent coordinates $(r, s, t)$, resulting in a matrix with nine rows and the same number of columns as in the matrix in (6). Next, to obtain the partial derivatives of lamina displacements with respect to the current position vector in the lamina CS, we reorder the rows to match $\begin{bmatrix}\frac{\partial u_1^l}{\partial r} & \frac{\partial u_1^l}{\partial s} & \frac{\partial u_1^l}{\partial t} & \frac{\partial u_2^l}{\partial r} & \frac{\partial u_2^l}{\partial s} & \frac{\partial u_2^l}{\partial t} & \frac{\partial u_3^l}{\partial r} & \frac{\partial u_3^l}{\partial s} & \frac{\partial u_3^l}{\partial t}\end{bmatrix}^T$, and then pre-multiply the sub-matrices by the inverse of the current Jacobian in the lamina CS (evaluated from $[J_a^l]^{-1} = [\,^{\tau}q_a][\,^{\tau}J_a]^{-1}$), to get (7). Note that:

1. the left super and subscripts ($\tau$ and $\beta$) are dropped for convenience,
2. $h_a$ is the nodal thickness (i.e. constant) thus its derivatives are zero,
3. the shape functions are 2D (i.e. independent of $t$), thus $\frac{\partial N_a}{\partial t} = 0$,
4. $\bar{t}$ is a constant and its derivatives are zero,
5. $\frac{\partial t}{\partial t} = 1$ and $\frac{\partial t}{\partial r} = \frac{\partial t}{\partial s} = 0$,
6. $\left.\begin{array}{l} g_{1i}^a = -\frac{1}{2}h_a\left(\vec{e}_1^{\,f}\right)_i^a \\ g_{2i}^a = -\frac{1}{2}h_a\left(\vec{e}_2^{\,f}\right)_i^a \end{array}\right\},\ i = 1,2,3$
7. $(t - \bar{t}) = \tilde{t}$.

$$\left\{\begin{array}{c}\left\{\begin{array}{c}\frac{\partial u_1^l}{\partial y_1^l}\\ \frac{\partial u_1^l}{\partial y_2^l}\\ \frac{\partial u_1^l}{\partial y_3^l}\end{array}\right\}\\ \left\{\begin{array}{c}\frac{\partial u_2^l}{\partial y_1^l}\\ \frac{\partial u_2^l}{\partial y_2^l}\\ \frac{\partial u_2^l}{\partial y_3^l}\end{array}\right\}\\ \left\{\begin{array}{c}\frac{\partial u_3^l}{\partial y_1^l}\\ \frac{\partial u_3^l}{\partial y_2^l}\\ \frac{\partial u_3^l}{\partial y_3^l}\end{array}\right\}\end{array}\right\} = \left[\begin{array}{c}
\cdots [J_a^l]^{-1}\begin{bmatrix}\frac{\partial N_a}{\partial r}(\vec{e}_1^{\,l})_1^a & \frac{\partial N_a}{\partial r}(\vec{e}_1^{\,l})_2^a & \frac{\partial N_a}{\partial r}(\vec{e}_1^{\,l})_3^a & \tilde{t}\left(\frac{-h_a}{2}\right)\frac{\partial N_a}{\partial r}\left((\vec{e}_1^{\,l})^a\cdot(\vec{e}_1^{\,f})^a\right) & \tilde{t}\left(\frac{-h_a}{2}\right)\frac{\partial N_a}{\partial r}\left((\vec{e}_1^{\,l})^a\cdot(\vec{e}_2^{\,f})^a\right)\\ \frac{\partial N_a}{\partial s}(\vec{e}_1^{\,l})_1^a & \frac{\partial N_a}{\partial s}(\vec{e}_1^{\,l})_2^a & \frac{\partial N_a}{\partial s}(\vec{e}_1^{\,l})_3^a & \tilde{t}\left(\frac{-h_a}{2}\right)\frac{\partial N_a}{\partial s}\left((\vec{e}_1^{\,l})^a\cdot(\vec{e}_1^{\,f})^a\right) & \tilde{t}\left(\frac{-h_a}{2}\right)\frac{\partial N_a}{\partial s}\left((\vec{e}_1^{\,l})^a\cdot(\vec{e}_2^{\,f})^a\right)\\ 0 & 0 & 0 & \left(\frac{-h_a}{2}\right)N_a\left((\vec{e}_1^{\,l})^a\cdot(\vec{e}_1^{\,f})^a\right) & \left(\frac{-h_a}{2}\right)N_a\left((\vec{e}_1^{\,l})^a\cdot(\vec{e}_2^{\,f})^a\right)\end{bmatrix}\cdots\\
\cdots [J_a^l]^{-1}\begin{bmatrix}\frac{\partial N_a}{\partial r}(\vec{e}_2^{\,l})_1^a & \frac{\partial N_a}{\partial r}(\vec{e}_2^{\,l})_2^a & \frac{\partial N_a}{\partial r}(\vec{e}_2^{\,l})_3^a & \tilde{t}\left(\frac{-h_a}{2}\right)\frac{\partial N_a}{\partial r}\left((\vec{e}_2^{\,l})^a\cdot(\vec{e}_1^{\,f})^a\right) & \tilde{t}\left(\frac{-h_a}{2}\right)\frac{\partial N_a}{\partial r}\left((\vec{e}_2^{\,l})^a\cdot(\vec{e}_2^{\,f})^a\right)\\ \frac{\partial N_a}{\partial s}(\vec{e}_2^{\,l})_1^a & \frac{\partial N_a}{\partial s}(\vec{e}_2^{\,l})_2^a & \frac{\partial N_a}{\partial s}(\vec{e}_2^{\,l})_3^a & \tilde{t}\left(\frac{-h_a}{2}\right)\frac{\partial N_a}{\partial s}\left((\vec{e}_2^{\,l})^a\cdot(\vec{e}_1^{\,f})^a\right) & \tilde{t}\left(\frac{-h_a}{2}\right)\frac{\partial N_a}{\partial s}\left((\vec{e}_2^{\,l})^a\cdot(\vec{e}_2^{\,f})^a\right)\\ 0 & 0 & 0 & \left(\frac{-h_a}{2}\right)N_a\left((\vec{e}_2^{\,l})^a\cdot(\vec{e}_1^{\,f})^a\right) & \left(\frac{-h_a}{2}\right)N_a\left((\vec{e}_2^{\,l})^a\cdot(\vec{e}_2^{\,f})^a\right)\end{bmatrix}\cdots\\
\cdots [J_a^l]^{-1}\begin{bmatrix}\frac{\partial N_a}{\partial r}(\vec{e}_3^{\,l})_1^a & \frac{\partial N_a}{\partial r}(\vec{e}_3^{\,l})_2^a & \frac{\partial N_a}{\partial r}(\vec{e}_3^{\,l})_3^a & \tilde{t}\left(\frac{-h_a}{2}\right)\frac{\partial N_a}{\partial r}\left((\vec{e}_3^{\,l})^a\cdot(\vec{e}_1^{\,f})^a\right) & \tilde{t}\left(\frac{-h_a}{2}\right)\frac{\partial N_a}{\partial r}\left((\vec{e}_3^{\,l})^a\cdot(\vec{e}_2^{\,f})^a\right)\\ \frac{\partial N_a}{\partial s}(\vec{e}_3^{\,l})_1^a & \frac{\partial N_a}{\partial s}(\vec{e}_3^{\,l})_2^a & \frac{\partial N_a}{\partial s}(\vec{e}_3^{\,l})_3^a & \tilde{t}\left(\frac{-h_a}{2}\right)\frac{\partial N_a}{\partial s}\left((\vec{e}_3^{\,l})^a\cdot(\vec{e}_1^{\,f})^a\right) & \tilde{t}\left(\frac{-h_a}{2}\right)\frac{\partial N_a}{\partial s}\left((\vec{e}_3^{\,l})^a\cdot(\vec{e}_2^{\,f})^a\right)\\ 0 & 0 & 0 & \left(\frac{-h_a}{2}\right)N_a\left((\vec{e}_3^{\,l})^a\cdot(\vec{e}_1^{\,f})^a\right) & \left(\frac{-h_a}{2}\right)N_a\left((\vec{e}_3^{\,l})^a\cdot(\vec{e}_2^{\,f})^a\right)\end{bmatrix}\cdots
\end{array}\right]\left\{\begin{array}{c}\vdots\\ u_1^a\\ u_2^a\\ u_3^a\\ \theta_1^a\\ \theta_2^a\\ \vdots\end{array}\right\} \tag{7}$$

Finally, for consistency with the rest of the document, the rows of (7) are manipulated such that the resulting $[^{\tau}_{\tau}B_L^l]$ produces the following vector of strains:

$$\{\varepsilon^l\} = \begin{bmatrix}\frac{\partial u_1^l}{\partial y_1^l} & \frac{\partial u_2^l}{\partial y_2^l} & \frac{\partial u_3^l}{\partial y_3^l} & \frac{\partial u_1^l}{\partial y_2^l} + \frac{\partial u_2^l}{\partial y_1^l} & \frac{\partial u_1^l}{\partial y_3^l} + \frac{\partial u_3^l}{\partial y_1^l} & \frac{\partial u_2^l}{\partial y_3^l} + \frac{\partial u_3^l}{\partial y_2^l}\end{bmatrix}^T \tag{8}$$

which, in Voigt notations is written as: $\{\varepsilon^l\} = [\varepsilon_1^l \quad \varepsilon_2^l \quad \varepsilon_3^l \quad \varepsilon_4^l \quad \varepsilon_5^l \quad \varepsilon_6^l]^T$.

## 2.5. *Constitutive relations in the lamina coordinate system for small strain analysis*

Considering a linear elastic material, Young's modulus and Poisson's ratio are employed to define the material tangent matrix (also called the constitutive tensor) in the lamina CS. In the UL formulation with small strains, the constitutive tensor relates the Cauchy stresses and Almansi strains in the lamina CS by [13,25]:

$$\{\,^{\tau}\sigma^l\} = [^{\tau}_{\tau}C^l]\{\,_{\tau}\varepsilon^l\} \tag{9}$$

As an adaptation, we evaluate the lamina Almansi strain tensor from $[\,_{\tau}\varepsilon^l] = \frac{1}{2}\left([I] - [^{\tau}_{\beta}F^l]^{-T}[^{\tau}_{\beta}F^l]^{-1}\right)$, where the left superscript $\tau$ denotes the current configuration, and the left subscript $\beta$ denotes the reference configuration. The deformation gradient of the current configuration with respect to the reference configuration in the lamina CS $[^{\tau}_{\beta}F^l]$ will be obtained from (4).

Dropping the left scripts for convenience, (9) is written using the Voigt indicial notations as: $\sigma_i^l = C_{ij}^l\varepsilon_j^l$. Note that in this formulation, index 3 corresponds to the direction normal to the surface of the shell. To apply the zero normal stress condition (i.e. $\sigma_3^l = 0$) on a linear elastic material, the third row of $[^{\tau}_{\tau}B_L^l]$ (corresponding to the normal strain component $\varepsilon_3^l$) and the third row and column of $[C^l]$ (corresponding to $\sigma_3^l$ and $\varepsilon_3^l$, respectively) must be removed. Consequently, $\varepsilon_3^l$ needs to be re-evaluated in terms of the remaining strain and stress components, such that [11,31]:

$$\varepsilon_3^l = \frac{-\left(\sum_{j=1,2,4,5,6} C_{3j}^l \varepsilon_j^l\right)}{C_{33}^l}. \tag{10}$$

The arrays with the third row and/or column removed are referred to as the "reduced" arrays, and are presented by a tilde over bar (e.g. $\{\tilde{\sigma}^l\}$). Substituting (10) in (9) and simplifying, gives the other components of the reduced Cauchy stress vector as:

$$\tilde{\sigma}_i^l = \sum_{j=1,2,4,5,6} \left( C_{ij}^l - \frac{C_{i3}^l C_{3j}^l}{C_{33}^l} \right) \varepsilon_j^l \tag{11}$$

Thus, writing (11) in matrix form, gives the reduced constitutive relation (satisfying the zero-normal stress condition) in the lamina CS, as:

$$\{ {}^{\tau}\tilde{\sigma}^l \} = [ {}^{\tau}_{\tau}\tilde{C}^l ] \{ {}_{\tau}\tilde{\varepsilon}^l \} \tag{12}$$

## 2.6. *Fiber length update algorithm for large membrane strains*

In the modified Mindlin-Reissner shell theory, the fiber inextensibility (constant fiber length) condition is invoked to avoid numerical ill-conditioned problems. Consequently, the fiber length parameters have not been included into the global nodal unknowns. However, this disqualifies the application of this theory in cases where large membrane strains may develop. To remove this limitation and update the fiber length parameters, the following procedure is used [2,21]. First, the strain tensor in the current configuration $[ {}_{\tau}\varepsilon^l ]$ obtained from the zero normal stress assumption must be transformed to the global CS using $[ {}_{\tau}\varepsilon ] = [ {}^{\tau}q ]^T [ {}_{\tau}\varepsilon^l ][ {}^{\tau}q ]$. Then, the mean value of $[ {}_{\tau}\varepsilon ]$ over the fiber is computed from:

$$[ {}_{\tau}\bar{\varepsilon} ] = \frac{1}{2}\left( \int_{-1}^{+1} [ {}_{\tau}\varepsilon ] \, dt \right). \tag{13}$$

Then, we project the mean component of strain (13) to produce the straining in the fiber direction. Recalling that $\hat{Y}$ denotes the unit vector in the fiber direction at the point in question, the transformation is done as ${}_{\tau}\bar{\varepsilon}^f = {}^{\tau}\hat{Y}^T [ {}_{\tau}\bar{\varepsilon} ] \, {}^{\tau}\hat{Y}$. Finally, we update the nodal fiber length parameters by ${}^{\tau}h_a \leftarrow {}^{\beta}h_a(1 + {}^{\tau}\bar{\varepsilon}_a^f)$. Note that the fiber length update is performed at the end of each iteration, which is after the element residual force and tangent array are completed. Therefore the update of $h_a$ lags one step behind the other kinematical quantities [21]. We also note that the thickness update results from the conservation of matter, and that volume evolution according to the constitutive relations is permitted.

## 2.7. *UL internal and external force vectors*

The internal force vector is obtained from the reduced Cauchy stress vector and the reduced strain-displacement transformation matrix, both within the lamina CS and in the current configuration. That is:

$$\{ {}^{\tau}_{\tau}f^{Internal} \} = \int_s \int_r \int_{-1}^{+1} [ {}^{\tau}_{\tau}\tilde{B}_L^l ]^T \{ {}^{\tau}\tilde{\sigma}^l \} J \, dt \, dr \, ds \, . \tag{14}$$

where, for accuracy in large deformations (nonlinear analysis), the Jacobian determinant $J$ is updated at each time step.

Due to only minor differences from the literature, the evaluation of the external force vectors is presented in Appendix D.

## 2.8. *UL mass matrix*

The mass matrix in the UL formulation is obtained from:

$$[ {}^{\tau}M^{consistent} ] = \int_s \int_r \int_{-1}^{+1} [ {}^{\tau}N ]^T [ {}^{\tau}N ] \, \rho \, J \, dt \, dr \, ds \tag{15}$$

where, for accuracy in large deformations, the Jacobian determinant $J$ and interpolation matrix $[N]$ are updated at each time step. For its convenience, the row sum technique can be used to lump the mass matrix:

$$ {}^{\tau}M_{ii} = \sum_{j=1}^{n} {}^{\tau}M_{ij}^{consistent} \tag{16} $$

## 2.9. *Critical time step*

An explicit time integration scheme is numerically stable as long as time step $\Delta t < \Delta t_{critical}$. The critical time step is equal to:

$$ \Delta t_{critical} = \frac{2}{\omega_n} \tag{17} $$

where $\omega_n$ is the maximum natural frequency in the finite element assemblage. This time step restriction applies to both linear and nonlinear systems [25]. The maximum natural frequency is expected to increase if the system stiffens and to decrease if the system softens. Natural frequencies are solutions of

$$ \det(-\omega^2[\,{}^{\tau}M_{ii}] + [{}^{\tau}_{\tau}K_L]) = 0, \tag{18} $$

where

$$ [{}^{\tau}_{\tau}K_L] = \int_s \int_r \int_{-1}^{+1} [{}^{\tau}_{\tau}\tilde{B}_L^l]^T [{}^{\tau}_{\tau}\tilde{C}^l][{}^{\tau}_{\tau}\tilde{B}_L^l]\, J\, dt\, dr\, ds\,. \tag{19} $$

The natural frequencies can be obtained by reducing (18) to a standard eigen value problem, and then using MATLAB's library routine to solve for the eigen vector $(\lambda)$ whose entries are the square of the natural frequencies $(\omega^2)_i$. That is: $\lambda = \mathrm{eig}([\,{}^{\tau}M_{ii}]^{-1}[{}^{\tau}_{\tau}K_L])$.

The advantages of using the eigen value approach to evaluate $\Delta t$ over the classic length over speed of sound method are threefold: 1) it results in larger $\Delta t$ values (see Tables 1 and 2), 2) it is applicable to elements with initially irregular and complex geometries, and 3) it is adaptable: $\Delta t$ needs to be evaluated for one element only (either the smallest element or the element undergoing the maximum deformation), and only as often as required by the application and the extent of the deformations; no adaptation is needed for small strain analysis in linear elasticity. Thus, the computational cost associated with this method is reasonable.

## 2.10. *Numerical integration*

We employed biquadratic lamina interpolation functions and the full (3×3) Gaussian rule for its simplicity and reliability [6,22]. This choice circumvented issues with reduced integration and selectively reduced integration, which are associated with zero-energy or kinematic modes (hourglassing) [26,29]. Hourglassing control, although possible, complicates computations, and stiffens the overall response of the system.

Considering fiber integrals, since the proposed CB shell element consists of one homogenous elastic layer, then the integrand is a smooth function of $t$ and the Gaussian quadrature is expected to be the most efficient. With the reference surface taken to be the midsurface $(\text{i.e. } \bar{t} = 0)$, we employed the 2-point Gauss rule to capture through-thickness variations [6,22].

As will be discussed in Section 4.1, full numerical integration (3×3×2) is one of the helpful factors in the prevention of shear and membrane locking in the present formulation.

## 2.11. *Operation count*

Considering that existing CB shell FEs were formulated in different programing languages and run on machines with dissimilar capabilities, the CPU times reported in the literature is not a neutral measure for efficiency comparisons. We believe the operation count (i.e. counting the total number of loops/iterations) to be a more equitable means of comparison. The FE formulations normally follow the main routine presented in Fig. 6, from which the number of operations can be evaluated [1–3,24,25,28,30]. Note that with the explicit method, there is no iteration per load increment.

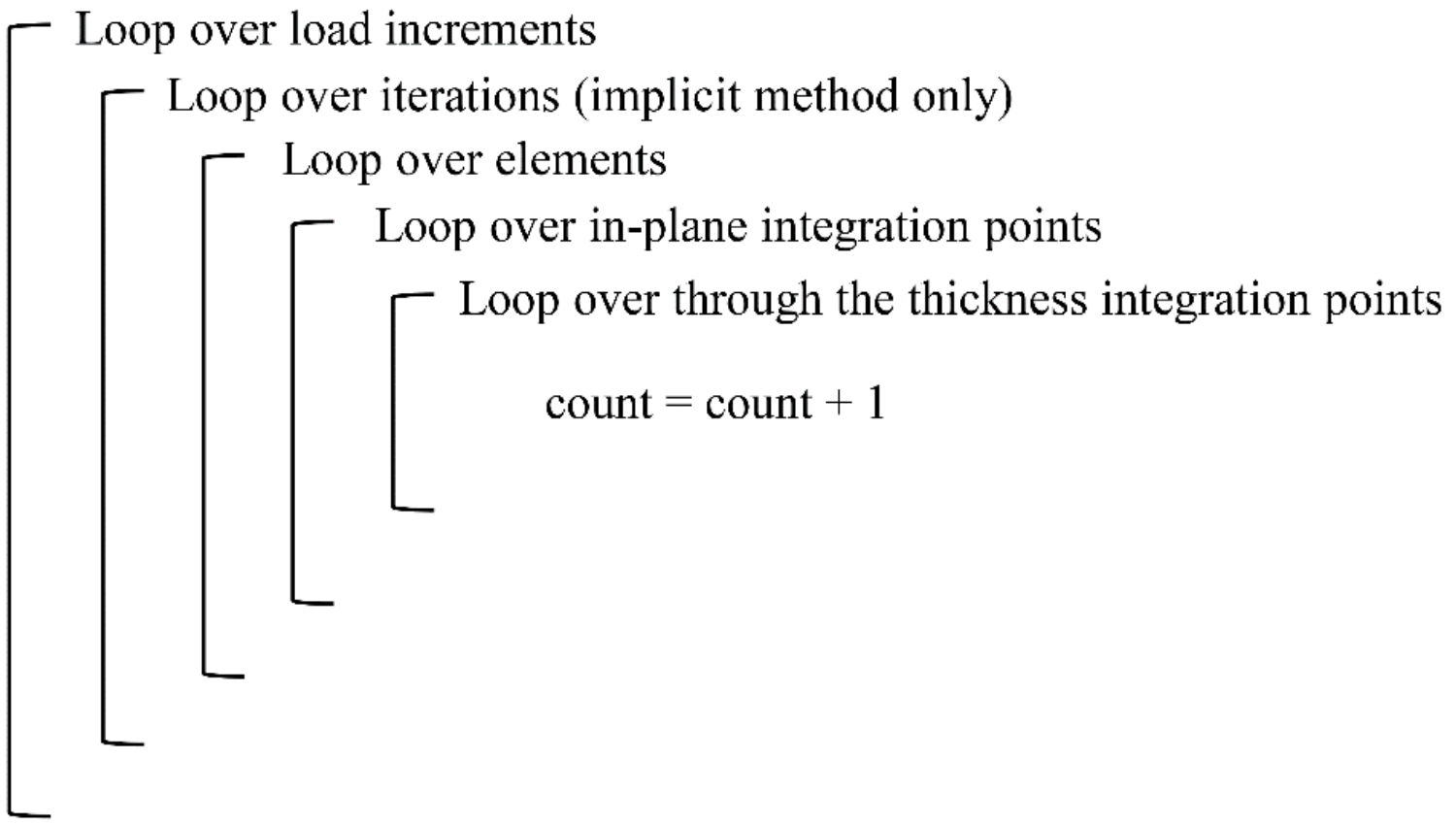

Fig. 6. Main routine used in FE formulations, from which, the number of operations is evaluated

## 3. Experiments

To establish the accuracy and efficiency of the present CB shell FE, it was submitted to a range of benchmark problems presented in many references.

### 3.1. *Experiment 1: Linear elastic, small bending deformations and rotations, small strains*

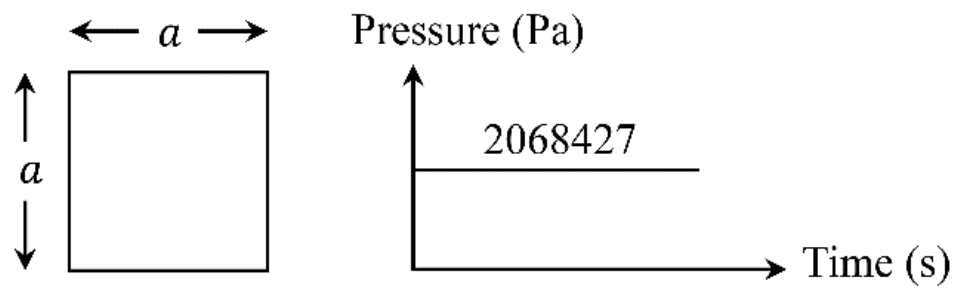

$h$ = Thickness = 0.0127 m
$a$ = Plate sides = 0.254 m
$E = 68.94 \times 10^9$ Pa
$v = 0.3$

Analytical linear static solution:

$$w_{max} = 0.00416\,\frac{12(1-v)Pa^4}{Eh^3} = 0.002794 \text{ m}$$

Fig. 7. Geometry, loading condition, and analytical linear static solution for Experiment 1

Table 1. Results of Experiment 1.

| | **Ref. [9]: 9 plates/quarter** | **Ref. [33]: 16 plates/quarter** | **Ref. [26]: 16 plates/quarter** | **Present: 1 shell/quarter** | **Present: 4 shells/quarter** |
|---|---|---|---|---|---|
| **Element type** | 8-noded | 4-noded | 4-noded | 9-noded | 9-noded |
| **Integration points** | 2×2×2 | NA | 1×1×1 | 3×3×2 | 3×3×2 |
| **Error (20)** | + 2.3 % | − 2.3 % | − 9.0 % | − 7.2 % | − 1.2 % |
| **Time integration** | Implicit | Implicit | Explicit | Explicit | Explicit |
| **Average Δt (s)** | $2.2\times10^{-5}$ | $2.2\times10^{-5}$ | $6.0\times10^{-6}$ | $8.4\times10^{-6}$ | $4.2\times10^{-6}$ |
| **Load increments** | 125 | 125 | 100 | 78 | 121 |
| **Total loops to max. deflection (Fig. 6)** | 9,000 | 2,000×NA | 1,600 | 1,404 | 8,712 |

A dynamic analysis of a simply supported aluminium square plate subjected to a step normal uniform pressure (Fig. 7) was carried out and the results were compared with the analytical linear static solution [8] and other finite elements introduced in [9,26,33] (Table 1). The maximum deflection in the center of the plate versus time (i.e. the dynamic response due to a step pressure) obtained from the present 1 shell and 4 shells per quarter meshes are presented in Fig. 8.

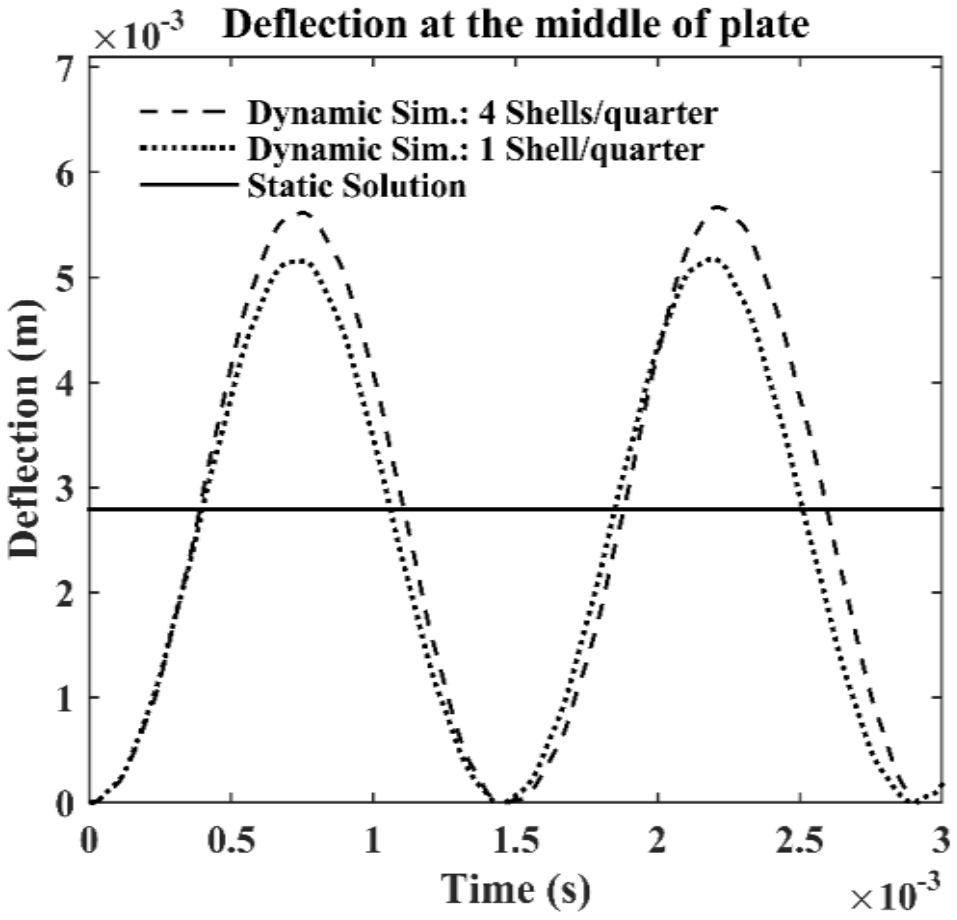

Fig. 8. Dynamic response of the simply supported plate due to step pressure (Experiment 1).

As expected from the step loading condition and absence of damping in the analysis, the dynamic response oscillated at constant amplitude with a magnitude almost equal to the analytical linear static solution (Fig. 8). Thus, the convergence error was calculated from:

$$error = \frac{amplitude - w_{max}^{analytical}}{w_{max}^{analytical}} \times 100\% \tag{20}$$

### 3.2. *Experiment 2: Elastic, large bending deformations and rotations, small strains*

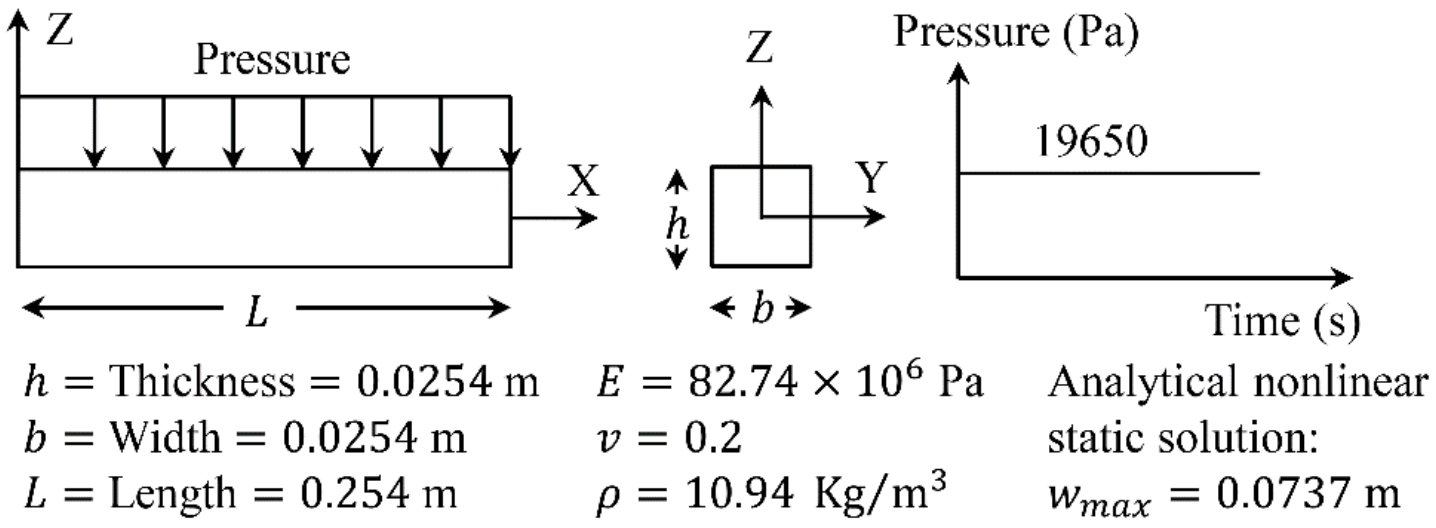

Fig. 9. Geometry, loading condition, and analytical nonlinear static solution for Experiment 2.

Table 2. Results of Experiment 2.

| ` | Ref. [29]: 5 plates | Present: 2 shells | Present: 3 shells |
|---|---|---|---|
| **Element type** | 8-noded | 9-noded | 9-noded |
| **Integration point** | 2×2×1 | 3×3×2 | 3×3×2 |
| **Error (20)** | + 3.4 % | − 4.3 % | − 2.8 % |
| **Time integration** | Explicit | Explicit | Explicit |
| **Average Δt (s)** | $2.0\times10^{-6}$ | $3.7\times10^{-6}$ | $3.6\times10^{-6}$ |
| **Total loops to max. deflection (Fig. 6)** | 33,000 | 25,056 | 39,042 |

A dynamic analysis of a cantilevered beam subjected to step normal uniform pressure (Fig. 9) was carried out. The large magnitude of the uniform pressure resulted in large bending deformations (nonlinear response), but small strains. The results obtained were compared with the analytical nonlinear static solution obtained from page 17 of [34] and another finite element introduced in [29] (Table 2).

### 3.3. *Experiment 3: Elastic, large, pure bending deformations and rotations, small strains*

The cantilevered beam shown in Fig. 10 was analyzed for its small strain, large displacement and large rotation (up to 90 degrees) response due to a concentrated end moment (i.e. large pure bending deformation), and the results were verified against the static analytical solution derived from page 54 of [34]. The idealization was done using regular and irregular 3- and 4- element (of the present 9-noded CB shell) meshes (Fig. 11). The normalized tip displacements (numerical over analytical) of the irregular 3-element mesh evaluated from this study and those tabulated in [17] are presented in Table 3. Fig. 12 shows the normalized axial, transverse and rotational displacements vs. moment parameters due to the applied tip moment for the 4 meshes considered.

Table 3. Results of Experiment 3.

| | **Ref. [17]: Irregular 3-element (4-noded) mesh** | | | **Present: Irregular 3-element (9-noded) mesh** | | |
|---|---|---|---|---|---|---|
| $\phi^{Analyt}$ | 18° | 45° | 72° | 18° | 45° | 72° |
| **Moment parameter** | 0.05 | 0.125 | 0.2 | 0.05 | 0.125 | 0.2 |
| $\phi^{Num}/\phi^{Analyt}$ | 0.95 | 0.84 | 0.76 | 1.05 | 1.02 | 0.93 |
| $U^{Num}/U^{Analyt}$ | 0.89 | 0.68 | 0.56 | 1.06 | 1.03 | 0.85 |
| $V^{Num}/V^{Analyt}$ | 0.95 | 0.86 | 0.81 | 1.02 | 1.01 | 0.94 |

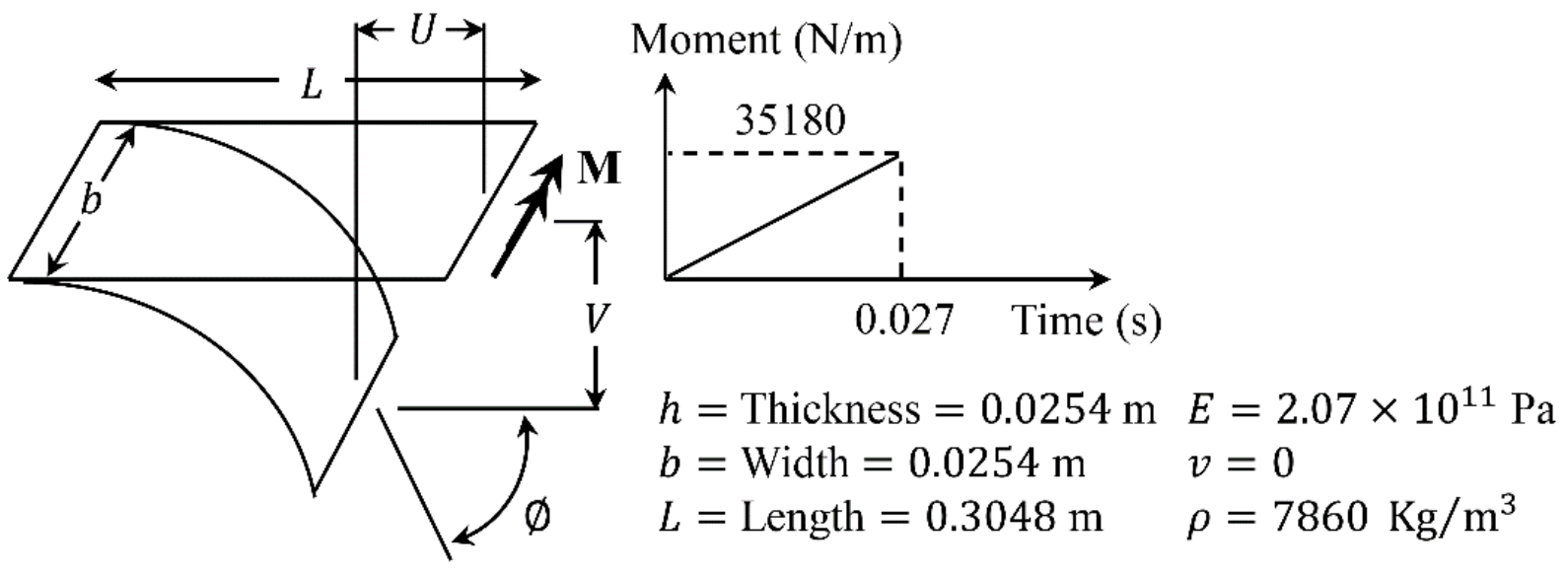

Fig. 10. Geometry, loading condition, and material properties for Experiment 3.

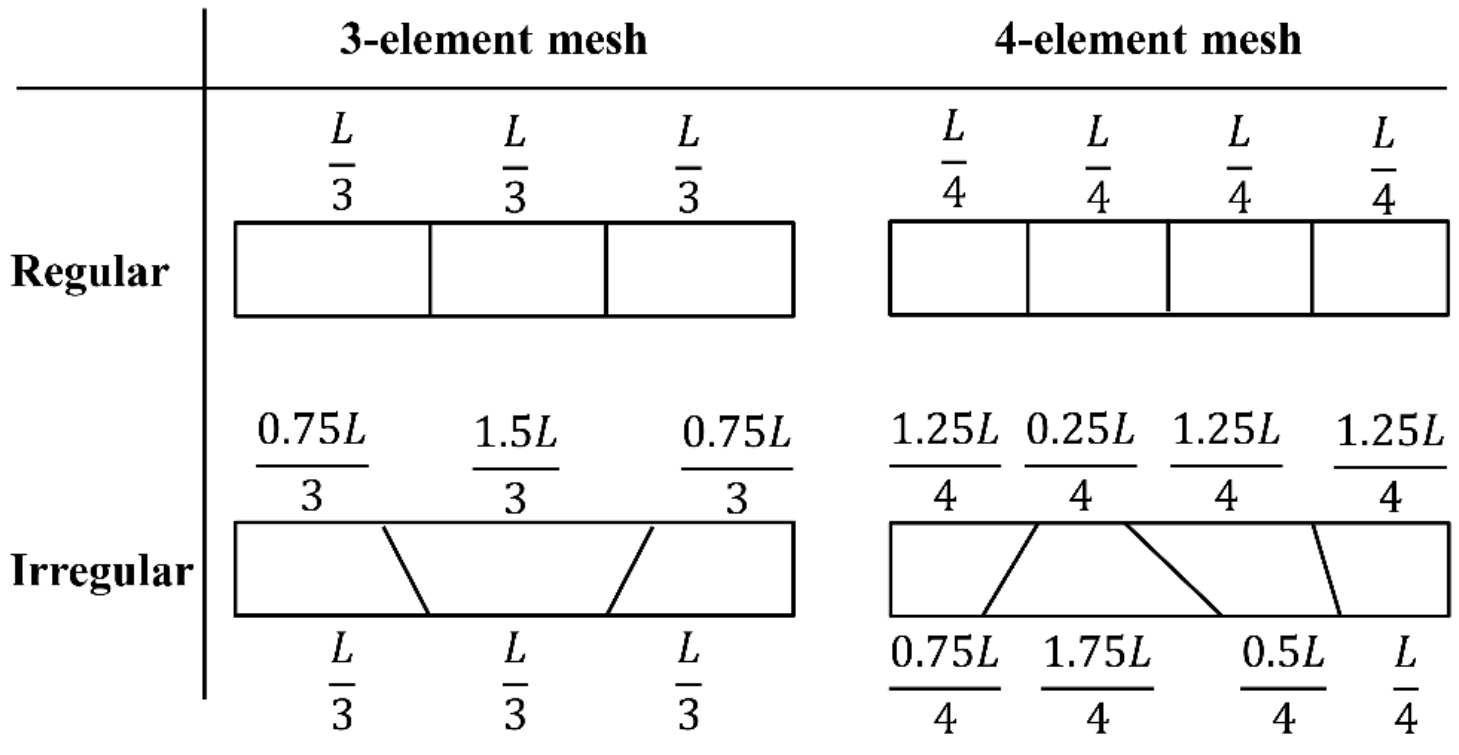

Fig. 11. Schematic of the four meshes considered in Experiment 3.

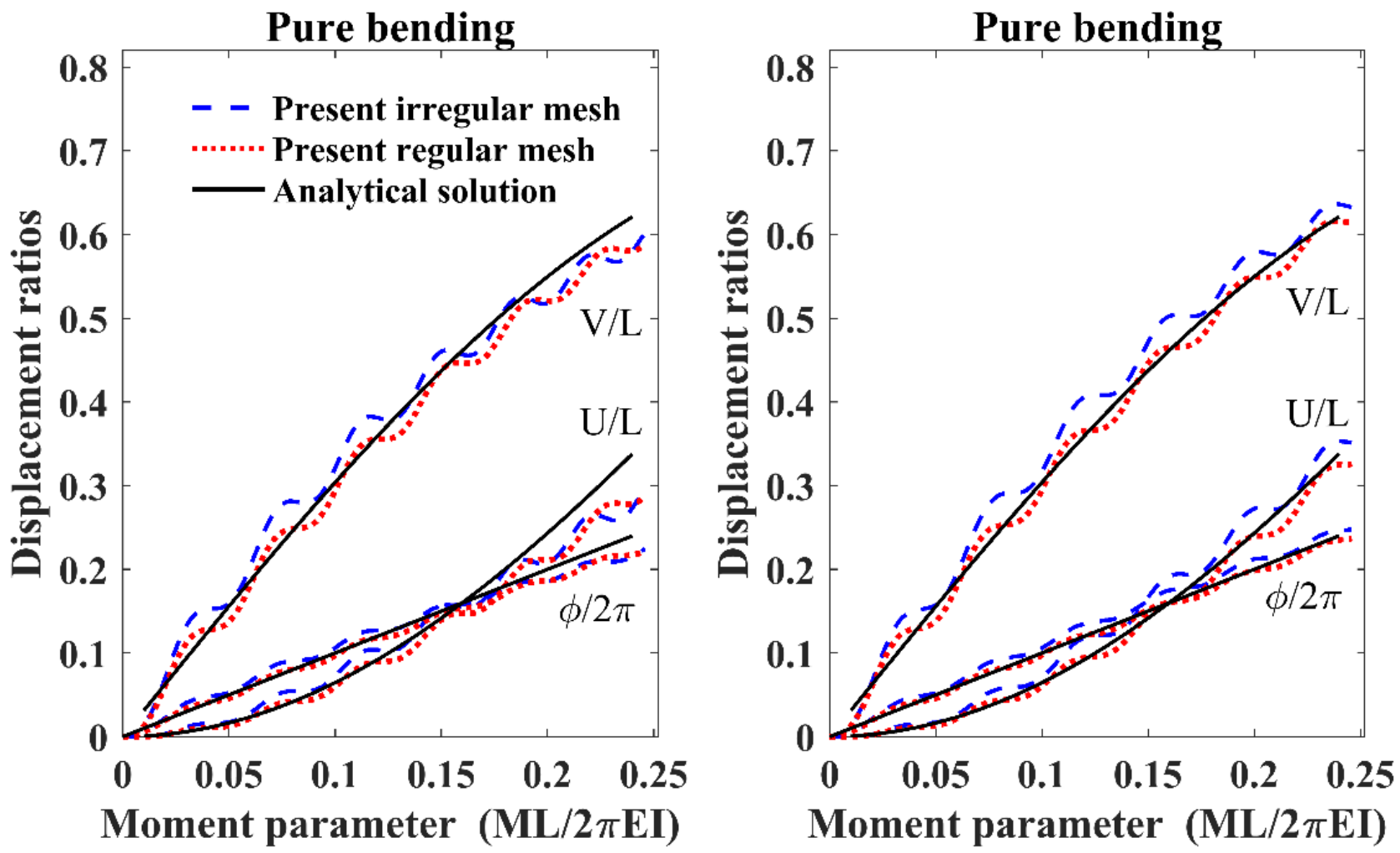

Fig. 12. Pure bending of a cantilever beam (Experiment 3). Left: 3-element mesh; Right: 4-element mesh.

### 3.4. *Experiment 4: Scordelis-Lo roof, initially singly-curved, membrane and bending deformation*

Fig. 13 illustrates the geometry, material properties, boundary conditions, and loading condition of a Scordelis-Lo roof (singly-curved shell structure) subjected to a uniform pressure in the vertical Z-direction. The theoretical value reported for the vertical displacement at the midpoint of the free edge is 0.3086 [35,36], but most elements converge to a slightly lower value (e.g. 0.302 in reference [37]). Convergence due to mesh refinement is presented in Fig. 14. The vertical displacement in region 1 with a 6×6 mesh is illustrated in Fig. 15.

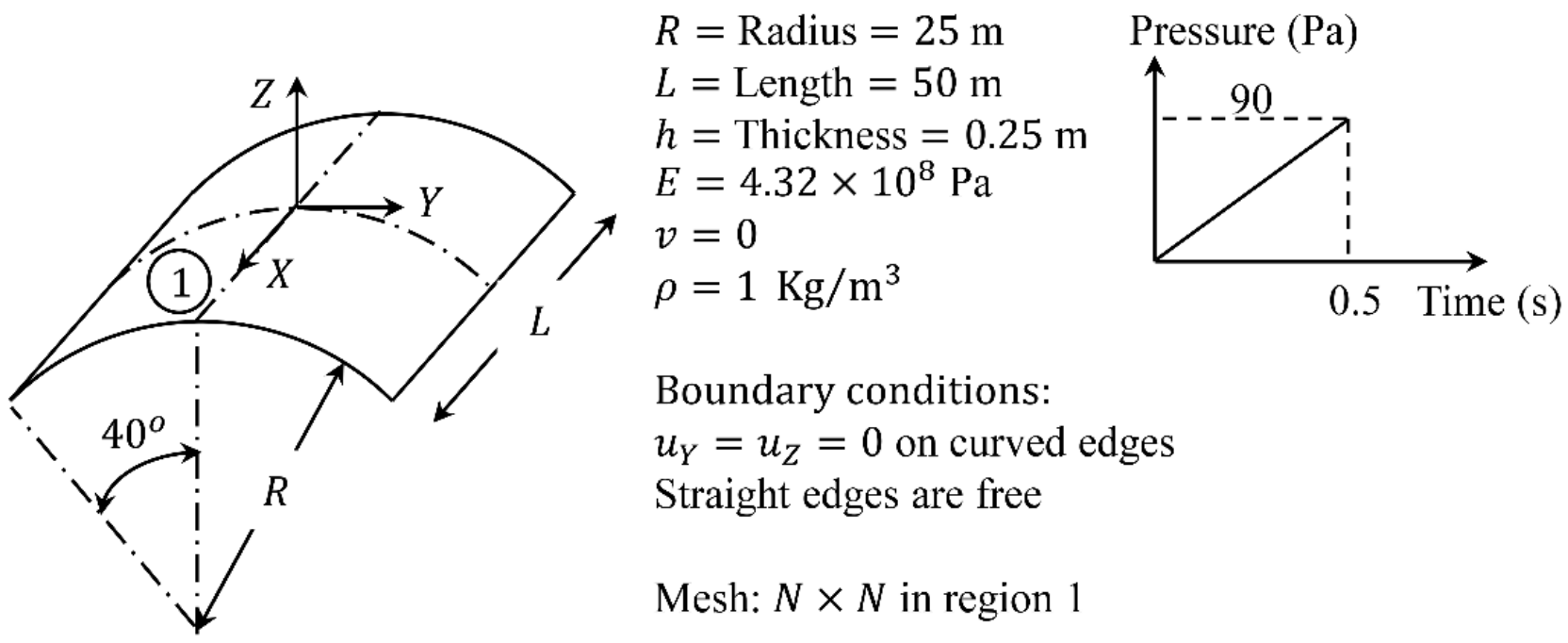

Fig. 13. Geometry, material properties, boundary conditions, and loading condition for a Scordelis-Lo roof (Experiment 4)

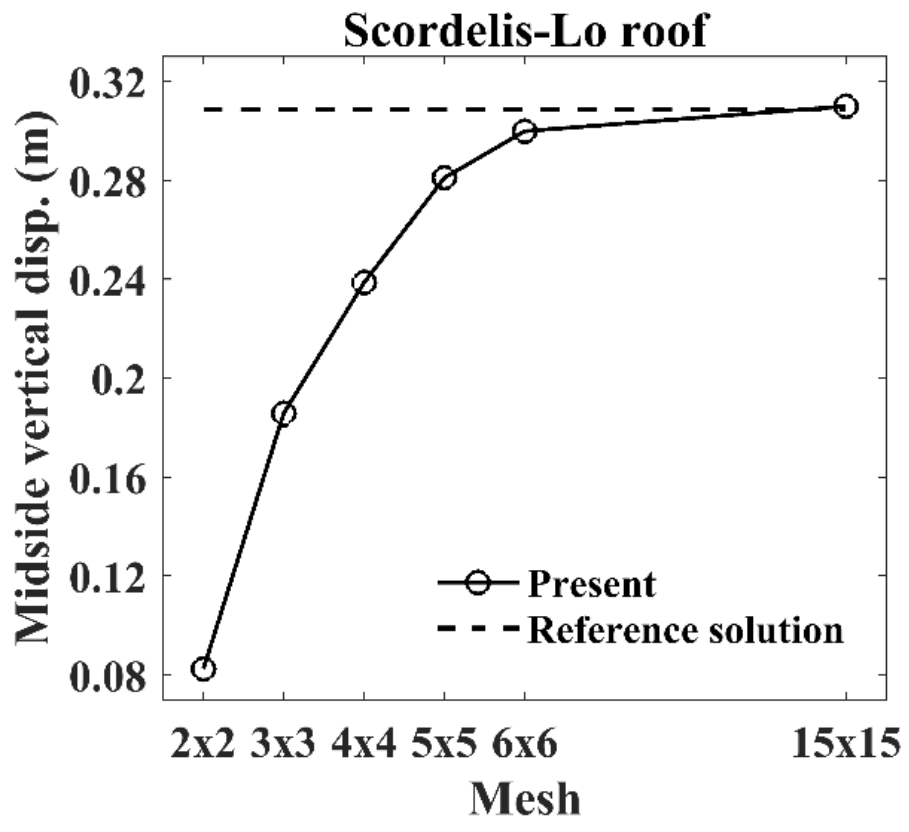

Fig. 14. Convergence of Experiment 4 to the analytical solution

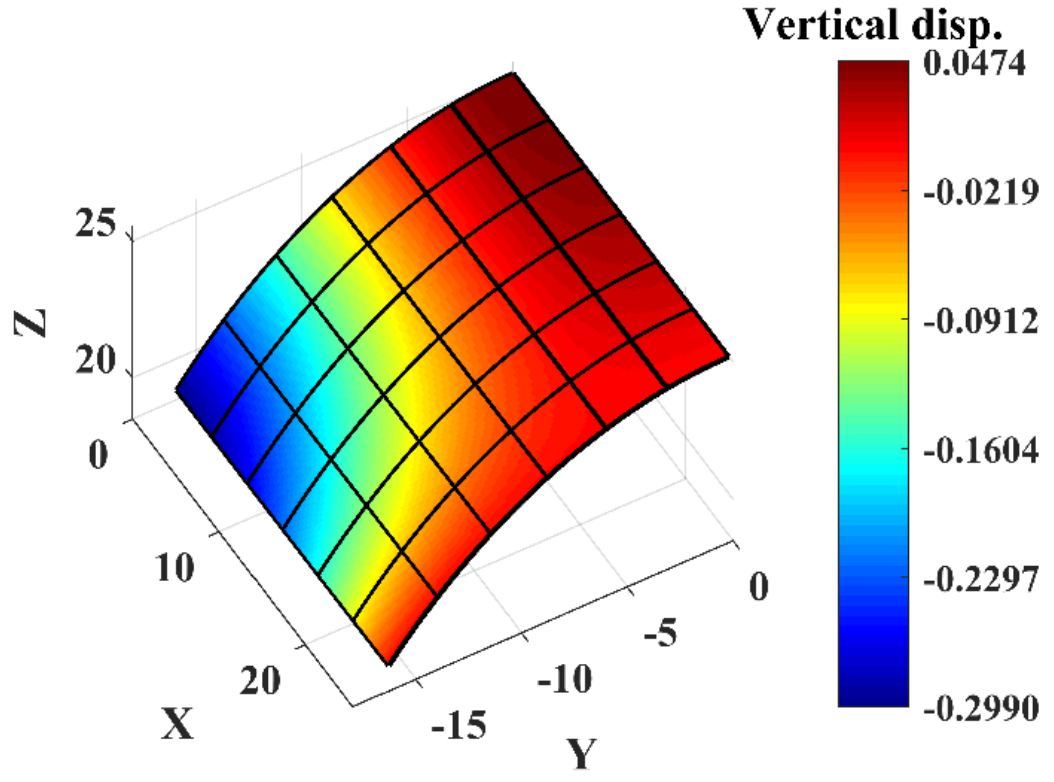

Fig. 15. Vertical displacement in region 1 of the Scordelis-Lo roof with a 6×6 mesh (Experiment 4). Units are in meters

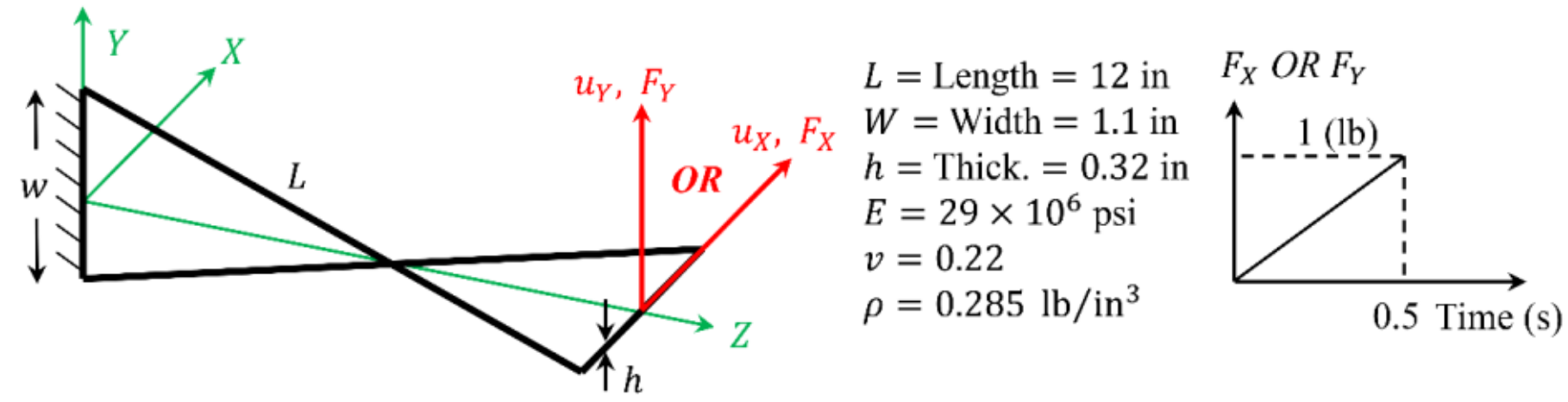

Fig. 16. Geometry, material properties, and loading condition of a cantilevered beam with an overall pre-twist of 90° (Experiment 5)

## 3.5. *Experiment 5: Large pre-twist, bending deformation in both planes*

The geometry, material properties, and loading condition of a cantilevered beam with an overall pre-twist of 90° are illustrated in Fig. 16. The theoretical solutions of tip deflections in the direction of loading are 0.005424 (in) and 0.001754 (in) for independent tip loading of one unit in the X- and Y-directions, respectively [35,38,39]. Convergence due to mesh refinement of the normalized displacements is presented in Fig. 17.

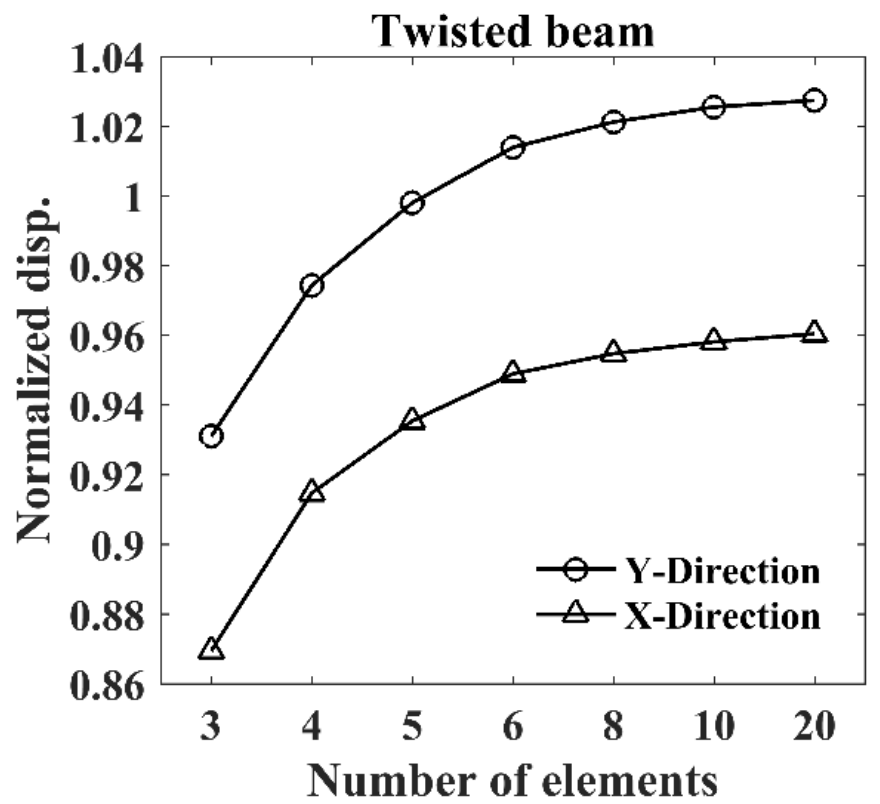

Fig. 17. Convergence of the normalized displacements in the X- and Y-directions due to mesh refinement (Experiment 5)

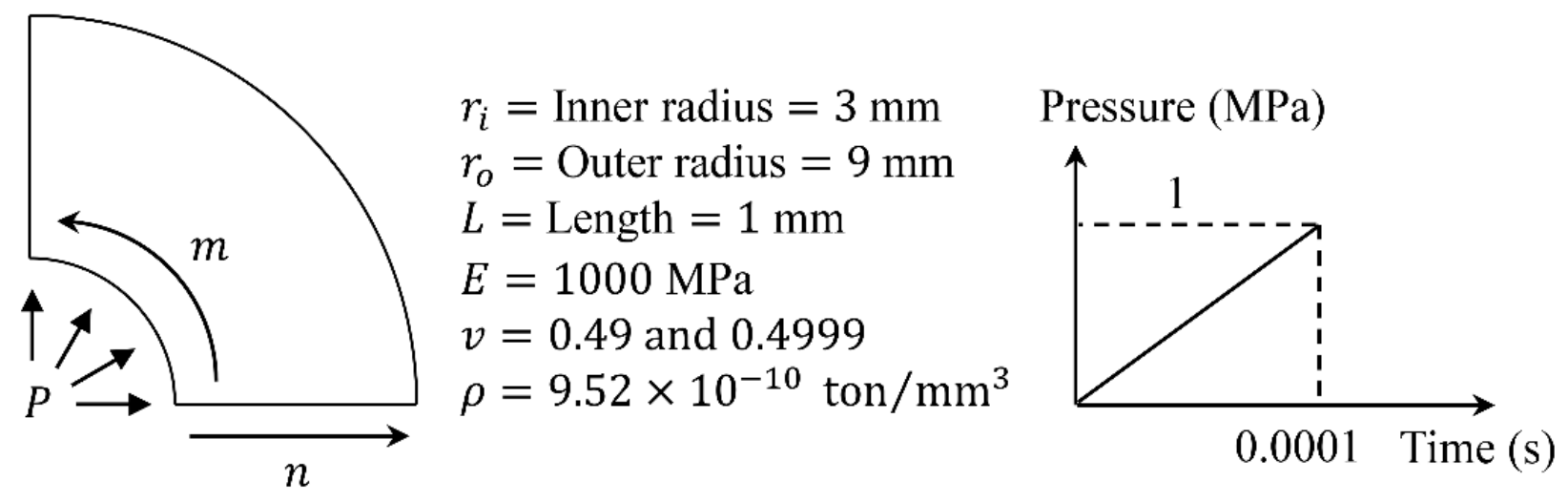

Fig. 18. Geometry, material properties, and loading condition for Experiment 6. $\boldsymbol{m}$ and $\boldsymbol{n}$ represent the number of elements in the circumferential and radial directions, respectively

Table 4. Results of Experiment 6: present errors obtained from (20), and computation costs for different mesh sizes and Poison's ratios. $r_i$ and $r_o$ , respectively, represent the internal and the external radii of the pipe.

| | **Present: 1×1** | | **Present: 2×2** | | **Present: 2×3** | |
|---|---|---|---|---|---|---|
| $v$ | **0.49** | **0.4999** | **0.49** | **0.4999** | **0.49** | **0.4999** |
| **Radial disp. error ($r_i$)** | 0.01% | 0.01% | 0.0% | 0.0% | 0.0% | 0.0% |
| **Circ. stress error ($r_i$)** | 16.5% | 16.9% | 4.4% | 4.6% | 1.5% | 1.5% |
| **Circ. stress error ($r_o$)** | 10.9% | 11.5% | 3.1% | 3.3% | 1.8% | 1.9% |
| **Rad. Stress error ($r_i$)** | −41.9% | −42.8% | −13.8% | −14.1% | −7.2% | −7.3% |
| **Average Δt (s)** | $4.3\times10^{-7}$ | $4.3\times10^{-8}$ | $1.8\times10^{-7}$ | $1.9\times10^{-8}$ | $1.4\times10^{-7}$ | $1.4\times10^{-8}$ |
| **Total loops to max. load (Fig. 6)** | $4.2\times10^{3}$ | $4\times10^{4}$ | $3.9\times10^{4}$ | $3.9\times10^{5}$ | $7.8\times10^{4}$ | $7.7\times10^{5}$ |

### 3.6. *Experiment 6: Thick-walled cylinder, linear elastic, small in-plane strains*

A wedge of an infinitely long thick-walled cylindrical pipe (Fig. 18, left), subjected to an internal pressure of 1 MPa (Fig. 18, right) is considered for the finite element analysis with an $m \times n$ mesh. An isotropic linear elastic with two values of Poisson's ratios, 0.49 and 0.4999, are separately considered to test the effect of near incompressibility. The plane-strain condition is assumed in the axial direction of the pipe which together with the radial symmetry confines the material in all but the radial direction. This problem was proposed in reference [35]. Table 4 details out the effect of mesh refinement and Poison's ratio on the percent errors. Fig. 19, illustrates the absence of stress jumps (in both circumferential and radial directions) across the elements for the 2×3 mesh.

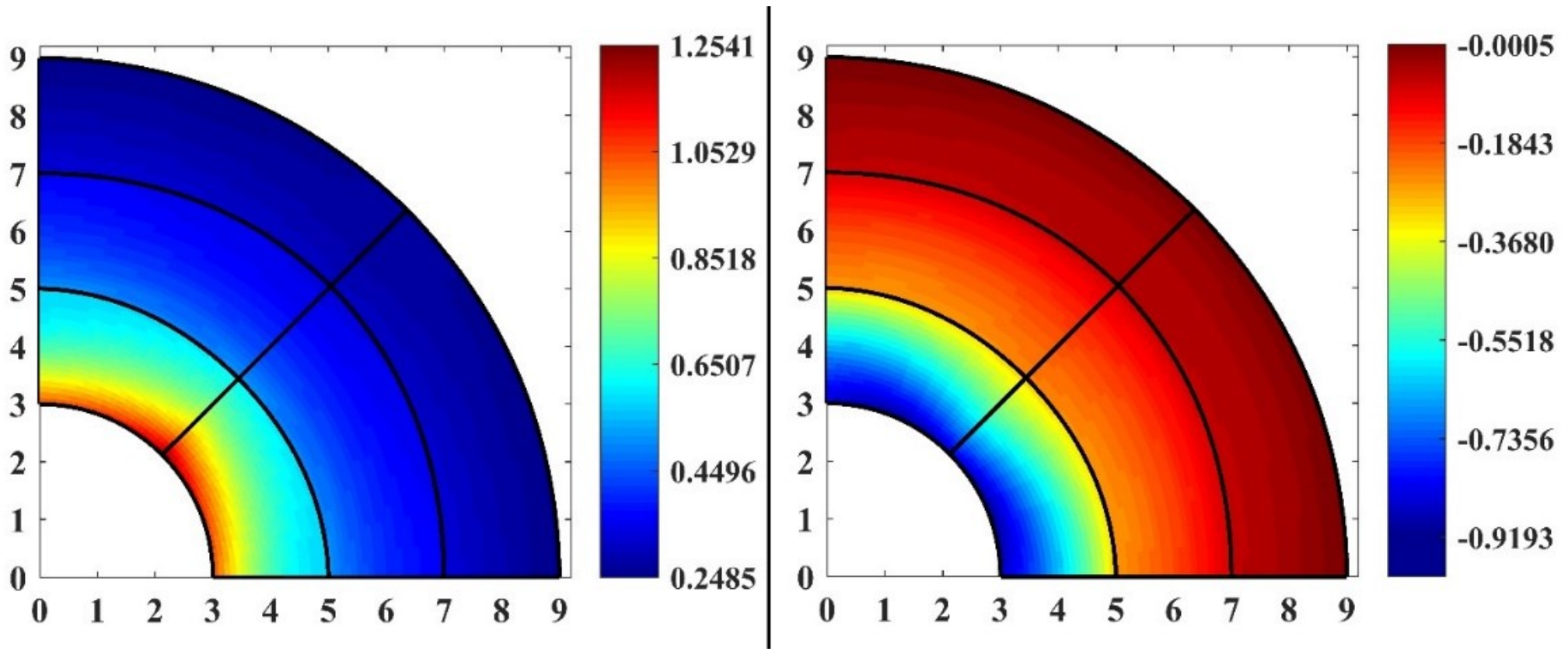

Fig. 19. Circumferential (left) and radial (right) stress distribution across the pipe wall under an internal pressure of 1 MPa, with = 0.49 , and for the 2×3 mesh (Experiment 6). Units for stress and radius are MPa and mm, respectively

## 4. Results and Discussion

### 4.1. *Benchmark experiments*

As presented in Table 1, use of only one of our 9-noded CB shell FE to model one quarter of the simply support plate under uniform normal pressure (Experiment 1), required the least total loops to reach the maximum deflection (i.e. computationally the most efficient) in comparison with [9,26,33], with an error (of −7.2%) more accurate than that of reference [26]. Furthermore, using the eigen value method (Section 2.9), time step $\Delta t$ was 38.2 times larger than those of reference [9] and reference [33], and 1.4 times faster than that of reference [26]. The references used the classical length over speed of sound equation for the critical time step.

By refining the mesh to four of the present 9-noded CB shell FE per quarter, we achieved an error of −1.2%. This was the most accurate compared to references [9,26,33], and slightly more efficient than reference [9]. The authors of [26] considered in-plane stresses and in-plane deformations only and used one-point quadrature rule for the numerical integration of their element. To avoid the resulting numerical ill-conditioning, they implemented hourglassing control. However, no information on the additional expense of this control on the computation time was provided. Overall, the best combination of accuracy and efficiency was achieved by employing four of the present CB shell FE per quarter of the plate.

To verify the accuracy and efficiency of our CB shell FE in large bending deformations (nonlinear analysis), we modeled a cantilevered beam subjected to a relatively large uniform pressure (Experiment 2). According to Table 2, using two of the present CB shell FE, $\Delta t$ was 1.9 times larger, the computation time was about 1.3 times shorter, and the error was 0.9% larger than those achieved in reference [29]. In the process of mesh refinement, we employed three of our CB shell FEs. This, in comparison with reference [29], resulted in an error 0.6% lower and a computation 1.2 times slower. Considering a good efficiency and accuracy combination, either the two-element or the three-element mesh of the present 9-noded CB shell FE did a good job of analysing the large deflection of this cantilever beam.

The high accuracy achieved in Experiments 1 and 2, even with few elements, was allowed by the fact that the limitation on the bending deformations/rotations was eliminated by the application of two independent CSs, and full numerical integration. Thus, the new CB shell FE has truly large bending deformation capabilities. Also, our calculated $\Delta t$ (Section 2.9) was many times larger than those obtained from dividing the length parameter by the sound speed. Although this approach is also used in [2] for solid (3D) shell elements, the present time step control is more efficient.

Through Experiment 2, we verified the accuracy and efficiency of the geometric nonlinear behaviour (i.e. large bending deformation) of the present CB shell FE when used as initially regular geometry. However, as reported in [8], higher order shell elements (i.e. quadrilateral elements that include more than four nodes) are generally sensitive to irregular geometries in the undeformed configuration, and lock in large (nonlinear) bending deformations. Hence, in Experiment 3, we studied the effect of using initially irregular elements for a cantilever beam under large deflections, and compared the results against the analytical solutions (reference [34], page 54). As presented in Fig. 12, left, the regular 3-element mesh yielded a very accurate response solution up to 45° , and as expected, we achieved a more accurate large deflection response (up to 90°) by refining the mesh to four regular elements (Fig. 12, right). Still, regardless of the number of elements used, the predictive capabilities of the present 9-noded CB shell FE were not affected by the irregularity of the mesh (Fig. 12). This claim is also supported by Table 3, in which the normalized (numerical/analytical) tip displacements (axial, transverse and rotational) obtained at different bending rotations and/or moment parameters from the present irregular 3-element mesh are very close to unity, and more accurate than those tabulated in reference [17]. Thus, not only was the formulation of the present 9-noded CB shell FE more straightforward than the mixed interpolation tensorial components used in [15–17], but it was also more accurate and more insensitive to initially irregular elements in large bending deformations. In addition, the accuracy of the tip (axial, transverse, and rotational) displacements verified that shear and membrane locking, which appear as false stiffening and false membrane extension in pure bending, are prevented. Assets in this regard are the evaluation of the lamina strain-displacement transformation matrix (Section 2.4), the implementation of constitutive relations in the lamina CS (Section 2.5), the derivation of the lamina CS according to Section 2.2 and Appendix B, and employing full numerical integration.

To further verify the performance of the present CB shell FE in handling curvatures, we conducted Experiment 4, in which, the Scordelis-Lo roof was subjected to a uniform vertical pressure. Both in-plane (membrane) and bending deformations contribute significantly to this singly-curved structure. Using a 6×6 mesh of the present 9-noded CB shell FE, the error in displacement was 2.8% (Fig. 14). This proves the present CB shell FE more accurate and more efficient than the 4-noded quadrilateral element presented in reference [37], and less accurate and less efficient than the 8-noded quadrilateral element in [35], where a mesh of 24×24, and a mesh of 2×2 were needed, respectively, to obtain the same error.

A common benchmark problem to determine the effects of warping in shell and solid finite elements is that of a tip-loaded, cantilevered beam, with 90° of overall pre-twist (Experiments 5) [35,37]. In the pre-twisted element, the lamina CS (the principal axes of the cross section) rotates along the element's length, and the pre-twist leads to a coupling of bending in both planes. Also, like in asymmetric bending, deflections of a pre-twisted cantilevered beam exhibit components both parallel and normal to the direction of loading. Comparison of the present FE results showed excellent convergence to the theoretical solution, indicating the accurate development of the stiffness matrix (19), bending ability in both planes, and insensitivity of the present CB shell FE to pre-twist. As illustrated in Fig. 17, displacement convergence, with errors of 2.1% and $-4.5\%$ in the Y- and X- direction, respectively, was achieved using 8 of the present CB shell FEs along the beam, where the warping of each element was 11.25°. Considering the 8-noded quadrilateral element presented in reference [35], convergence with an error of 2% in both directions was accomplished using 12 elements (7.5° warping per element) along the beam.

In Experiment 6, a 90° section of a thick-walled cylinder was chosen to test accuracy of the present CB shell FE, and its insensitivity to volumetric locking in 2D plane strain as $\upsilon \rightarrow 0.5$ (nearly incompressible materials). According to Table 4, radial displacement convergence was achieved with a 1×1 mesh, proving the present CB shell FE more accurate than the 8-noded quadrilateral element presented in [35], where a minimum of 1×5 mesh for a 10° section is needed to achieve the same result. Despite the displacement convergence in a 1×1 mesh of the present CB shell FE, the errors on stresses in the circumferential and radial directions were still large. Thus, to reduce the error on the radial stresses at the inner surface (the largest error) from $-41.9\%$ to $-7.2\%$, we refined the mesh to 2×3 elements, which made the computation only 2 times more expensive (Table 4). In reference [35], any error less than, or equal to, 10% is ranked as "very good" for this experiment. Thus, the best accuracy and efficiency combination was attained using a 2×3 mesh. The influence of increasing Poison's ratio to 0.4999 on the errors presented was negligible (a maximum of 0.2% on average), it did not cause volumetric locking, but decreased $\Delta t$ by one-tenth. Finally, the absence of jumps in the distribution of the radial and circumferential stresses across the elements (another measure of convergence) was verified for the 2×3 mesh (Fig. 19).

### 4.2. *Contributions of the new CB shell FE*

Compared to existing CB shell FE elements, the most important contributions of the new CB shell FE lie in its accuracy and efficiency in analyzing bending deformations in two planes (Experiments 1), large bending deformation (Experiment 2), large pure bending deformation of initially regular and irregular (distorted) elements (Experiment 3), combined membrane and bending deformations of initially curved geometries (Experiment 4), and combined bending deformation of pre-twist cantilever beam (Experiments 5). These characteristics were also obtained using a coarser mesh (i.e. fewer, larger elements). In essence, the improved accuracy using a coarser mesh reflect that shear and membrane locking, which respectively, appear as false stiffening and false membrane extension in pure bending, are avoided in the formulation of the present CB shell FE. In addition, we verified that the new CB shell FE prevents incompressible (volumetric) locking in 2D plane strain as the Poisson's ratio $\upsilon \rightarrow 0.5$ (Experiment 6). The new CB shell FE achieves good accuracy without the need for equilibrium iterations as in implicit methods, while reducing the number of operations and increasing the time step ($\Delta t$) (Experiments 1 and 2). Provided that all the CB shell elements (considered for comparison in the experiments) are implemented using the same programming language and the tests run on the same machine, both factors are indicators of reduced computational time.

### 4.3. *Limitations*

The new thick CB shell FE was implemented in the MATLAB environment due to its availability and convenience, but the element could be made numerically more efficient by using a faster programming language. Considering that two through the thickness integration points are necessary for the prevention of shear locking in bending deformations, the efficiency of the element could also be doubled by reducing the number of through the thickness integration points from 2 to 1 when only in-plane deformations are modelled (e.g. Experiment 6). Moreover, to further leverage the small number of elements required per model, the loop over the elements (Fig. 6) could be circumvented by running the elemental calculations in parallel, thereby reducing the number of operations and making the code even more efficient. However, the successful implementation of the new element with anisotropic, nearly incompressible, hyperelastic material properties to represent soft biological tissues still needs to be demonstrated; for the sake of space and of validating the foundations first, it was not included here. In a companion article (Part 2), constitutive relations needed for modelling rubber-like materials and soft biological tissues are presented, and applications of the present thick CB shell FE in the biomedical realm are further verified.

## 5. Conclusion

Geometrically and materially nonlinear FE analysis requires a CB shell FE that is accurate, reliable and versatile. The benchmark experiments performed proved the present thick CB shell FE to be powerful in handling combined large in-plane (membrane) and large bending deformations, initially irregular (distorted) elements and curved geometries, allowing for relatively coarse meshes even for complex geometries, and achieving comparatively fast computational times. The new element was also insensitive to three types of locking (shear, membrane and volumetric), and pre-twist (i.e. warping effects). Overall, from its robust derivation and application of two independent CSs, proper development of the stiffness and strain-displacement transformation matrices and constitutive relations in the lamina CS, accurate evaluation of the thickness changes, as well as optimal critical time step calculation, the new thick CB shell FE appeared as a promising addition to the list of existing CB shell FEs.

## Acknowledgements

This work was supported by the Natural Sciences and Engineering Research Council of Canada for Discovery Grant 312065-2012 (M.R.L).

## Appendix A: Derivation of the displacement of a director at a specific node $\hat{U}_a$

The fibers are assumed to be inextensible (in kinematics only), meaning that they can rotate but cannot stretch or contract. Therefore:

$$\left\|\hat{Y}_a\right\| = 1 \qquad (A.1)$$

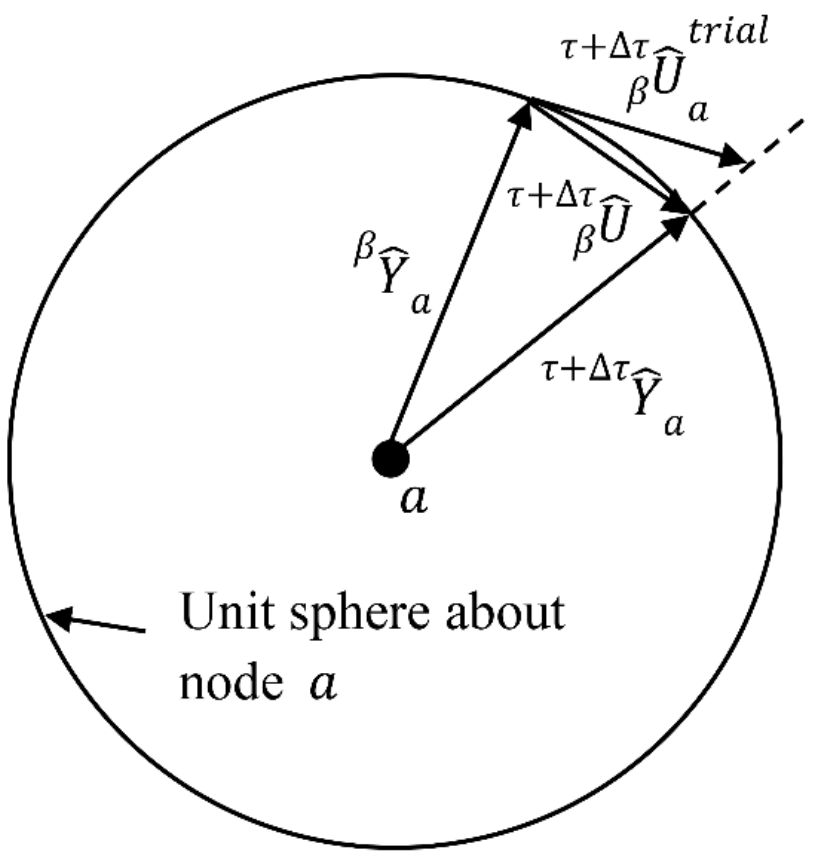

Fig. A.1. Nodal fiber inextensibility condition maintained by radial return normalization

A trial value of $\hat{U}_a$ (i.e $\hat{U}_a^{trial}$) is calculated from (A.4) and then projected radially (Fig. A.1) to maintain (A.1). The steps are as follows:

$$ {}^{\tau+\Delta\tau}\hat{Y}_a = \frac{\left({}^{\beta}\hat{Y}_a + {}^{\tau+\Delta\tau}_{\beta}\hat{U}_a^{trial}\right)}{\left\| {}^{\beta}\hat{Y}_a + {}^{\tau+\Delta\tau}_{\beta}\hat{U}_a^{trial} \right\|} \quad (A.2) $$

$$ {}^{\tau+\Delta\tau}_{\beta}\hat{U}_a = {}^{\tau+\Delta\tau}\hat{Y}_a - {}^{\beta}\hat{Y}_a \quad (A.3) $$

The trial value of the displacement of a director at a specific node, ${}^{\tau+\Delta\tau}_{\beta}\hat{U}_a^{trial}$, is obtained from:

$$ \begin{Bmatrix} \left({}^{\tau+\Delta\tau}_{\beta}\hat{U}_1^{trial}\right)_a \\ \left({}^{\tau+\Delta\tau}_{\beta}\hat{U}_2^{trial}\right)_a \\ \left({}^{\tau+\Delta\tau}_{\beta}\hat{U}_3^{trial}\right)_a \end{Bmatrix} = \begin{bmatrix} -\left({}^{\beta}\vec{e}_1^{\,f}\right)_1^a & -\left({}^{\beta}\vec{e}_2^{\,f}\right)_1^a \\ -\left({}^{\beta}\vec{e}_1^{\,f}\right)_2^a & -\left({}^{\beta}\vec{e}_2^{\,f}\right)_2^a \\ -\left({}^{\beta}\vec{e}_1^{\,f}\right)_3^a & -\left({}^{\beta}\vec{e}_2^{\,f}\right)_3^a \end{bmatrix} \begin{Bmatrix} {}^{\tau+\Delta\tau}_{\beta}\theta_1^a \\ {}^{\tau+\Delta\tau}_{\beta}\theta_2^a \end{Bmatrix}. \quad (A.4) $$

**Appendix B: Derivation of the lamina CS**

The bases of the lamina CS are obtained through the following procedure:

$$ \vec{e}_r = \frac{y_{,r}}{\|y_{,r}\|} \quad \text{and} \quad \vec{e}_s = \frac{y_{,s}}{\|y_{,s}\|}, \quad (B.1) $$

where $y$ denotes the current position vector in the global CS and the comma represents partial differentiation. Then,

$$ \vec{e}_3^{\,l} = \frac{\vec{e}_r \times \vec{e}_s}{\|\vec{e}_r \times \vec{e}_s\|}. \quad (B.2) $$

Following [31], the vectors tangent to the lamina are selected such that the angle between $\vec{e}_1^{\,l}$ and $\vec{e}_r$ (i.e. the vector tangent to $r$) is equal to the angle between $\vec{e}_s$ and $\vec{e}_2^{\,l}$. In addition, the $\vec{e}_1^{\,l}$ and $\vec{e}_2^{\,l}$ basis is as close as possible to the $\vec{e}_r$ and $\vec{e}_s$ basis.

$$ \vec{e}_1^{\,l} = \frac{\sqrt{2}}{2}\left(\vec{e}_\alpha - \vec{e}_\beta\right), \quad (B.3) $$

$$ \vec{e}_2^{\,l} = \frac{\sqrt{2}}{2}\left(\vec{e}_\alpha + \vec{e}_\beta\right), \quad (B.4) $$

where $\vec{e}_\alpha$ and $\vec{e}_\beta$ are obtained from:

$$ \vec{e}_\alpha = \frac{0.5 \times \left(\vec{e}_r + \vec{e}_s\right)}{\|0.5 \times \left(\vec{e}_r + \vec{e}_s\right)\|} \quad \text{and} \quad \vec{e}_\beta = \frac{\vec{e}_3^{\,l} \times \vec{e}_\alpha}{\|\vec{e}_3^{\,l} \times \vec{e}_\alpha\|}. \quad (B.5) $$

Note that the unit vectors of the lamina CS are functions of the parent CS and that they can be obtained for any specified $r, s, t$ value.

Note that, for simplicity in [7–23,25], $\vec{e}_1^{\,l}$ is bound with the tangent of the parent curvilinear coordinate $r$ (i.e. $\vec{e}_1^{\,l} = \vec{e}_r$ ). Thus, the evaluation of the lamina CS in that approach is very inaccurate when large distortions appear or when initially irregular elements are used.

**Appendix C: Derivation for the fiber CS**

As shown in Fig. 4 $\hat{Y}_i$ (where $i = 1, 2, 3$) denotes the projections of the current director $(\hat{Y})$ on the global Cartesian coordinate basis. Therefore, $\hat{Y} = |\hat{Y}_1|\vec{e}_1 + |\hat{Y}_2|\vec{e}_2 + |\hat{Y}_3|\vec{e}_3$.

So far, the direction of $\vec{e}_3^{\,f}$ has been determined to be coincident with the fiber direction $\left(\vec{e}_3^{\,f} = \hat{Y}\right)$. In order to determine the direction of $\vec{e}_1^{\,f}$ and $\vec{e}_2^{\,f}$, authors in [22] suggested the following algorithm:

*Step 1:* let $b_i = |\hat{Y}_i|,\ i = 1, 2, 3$
*Step 2:* $j = 1$
*Step 3:* If $b_1 > b_3$, then $b_3 = b_1$ and $j = 2$
*Step 4:* If $b_2 > b_3$, then $j = 3$
*Step 5:* $\vec{e}_3^{\,f} = \hat{Y}$
*Step 6:* $\vec{e}_2^{\,f} = \left(\hat{Y} \times \vec{e}_j\right) / \|\hat{Y} \times \vec{e}_j\|$
*Step 7:* $\vec{e}_1^{\,f} = \left(\vec{e}_2^{\,f} \times \hat{Y}\right)$

The above algorithm ensures that the obtained orthonormal fiber basis $\left(\vec{e}_1^f, \vec{e}_2^f, \vec{e}_3^f\right)$ satisfies the condition that if $\hat{Y}$ is close to $\vec{e}_3$ , and then $\vec{e}_1^f, \vec{e}_2^f, \vec{e}_3^f$ will be close to $\vec{e}_1, \vec{e}_2, \vec{e}_3$ respectively.

**Appendix D: Evaluation of external force vectors**

Two of the most common external force vectors are the body and the surface forces.

The body force vector is:

$$\left\{{}_{\beta}^{\tau}R^{B}\right\} = \int_{s}\int_{r}\int_{-1}^{+1}[\,{}^{\tau}N]^{T}\left\{{}_{\beta}^{\tau}f^{B}\right\}\rho\, J\, dt\, dr\, ds \qquad \text{(D.1)}$$

where, $\left\{{}_{\beta}^{\tau}f^{B}\right\}$ is the body force vector (per unit mass), $\beta$ denotes the reference configuration. In addition, $J$ and $[\,{}^{\tau}N]$ must be updated for each configuration for better accuracy, especially when large deformations and large rotations are considered. This fact is neglected in many of the references.

The surface force vector is:

$$\left\{{}_{\beta}^{\tau}R^{S_f}\right\} = \int_{s}\int_{r}[\,{}^{\tau}N]^{T}\left\{{}_{\beta}^{\tau}f^{S_f}\right\}J_s\, dr\, ds\,, \quad t = \begin{cases}+1 & \text{top}\\ -1 & \text{bottom}\end{cases} \qquad \text{(D.2)}$$

where, $\left\{{}_{\beta}^{\tau}f^{S_f}\right\}$ is the surface force vector (per unit surface area), and the current surface Jacobian is obtained from $J_s = \left\|\begin{pmatrix}\frac{\partial y_1}{\partial r} & \frac{\partial y_2}{\partial r} & \frac{\partial y_3}{\partial r}\end{pmatrix} \times \begin{pmatrix}\frac{\partial y_1}{\partial s} & \frac{\partial y_2}{\partial s} & \frac{\partial y_3}{\partial s}\end{pmatrix}\right\|$.

**References**

1. *ABAQUS User: Shell Elements in ABAQUS/Explicit*. SIMULIA; 2005.
2. *LS-DYNA Theoretical*. Livermore software Technology Corporation; 2011.
3. Segal G. *Sepra Analysis: Programmers Guide*. Netherlands: Sepran; 2010.
4. Einstein DR, Reinhall P, Nicosia M, Cochran RP, Kunzelman K. Dynamic finite element implementation of nonlinear, anisotropic hyperelastic biological membranes. *Comput Methods Biomech Biomed Engin*. 2003;6(March 2016):33-44. doi:10.1080/1025584021000048983.
5. Weinberg EJ, Kaazempur Mofrad MR. A finite shell element for heart mitral valve leaflet mechanics, with large deformations and 3D constitutive material model. *J Biomech*. 2007;40(3):705-711. doi:10.1016/j.jbiomech.2006.01.003.
6. Belytschko T, Liu WK, Moran B. Nonlinear Finite Elements for Continua and Structures. 2000.
7. Ahmad S, Irons BM, Zienkiewicz OC. Analysis of thick and thin shell structures by curved finite elements. *Int J Numer Methods Eng*. 1970;2(3):419-451. doi:10.1002/nme.1620020310.
8. Bathe K-J, Dvorkin E, Ho LW. Our discrete-Kirchhoff and isoparametric shell elements for nonlinear analysis-An assessment. *Comput Struct*. 1983;16(1-4):89-98. doi:10.1016/0045-7949(83)90150-5.
9. Bathe K jürgen, Bolourchi S. A geometric and material nonlinear plate and shell element. *Comput Struct*. 1980;11(1-2):23-48. doi:10.1016/0045-7949(80)90144-3.
10. Kanok-nukulchai W. A simple and efficient finite element for general shell analysis. *Int J Numer …*. 1979;14(April 1978):179-200. doi:10.1002/nme.1620140204.
11. Kiendl J, Hsu MC, Wu MCH, Reali A. Isogeometric Kirchhoff-Love shell formulations for general hyperelastic materials. *Comput Methods Appl Mech Eng*. 2015;291:280-303. doi:10.1016/j.cma.2015.03.010.
12. Krakeland B. *Large Displacement Analysis of Shells Considering Elastic-Plastic and Elasto-Viscoplastic Materials*. Norway: The Norwegian Institute of Technology; 1977.
13. Bathe K-J, Ramm E, Wilson E. Finite Element Formulations for Large Deformation Dynamic Analysus. *Int J Numer Methods Eng*. 1975;9:353-386. doi:10.1002/nme.1620090207.
14. Ramm E. A plate/shell element for large deflections and rotations. *Formul Comput algorithms finite Elem Anal*. 1977:264-293.
15. Dvorkin EN, Bathe K-J. A continuum mechanics based four-node shell element for general non-linear analysis. *Eng Comput*. 1984;1(1):77-88. doi:10.1108/eb023562.
16. Bucalem ML, Bathe K-J. Higher-order MITC general shell elements. *Int J Numer Methods Eng*. 1993;36(March):3729-3754. doi:10.1002/nme.1620362109.
17. Dvorkin EN. Nonlinear analysis of shells using the MITC formulation. *Arch Comput Methods Eng*. 1995;2(2):1-50. doi:10.1007/BF02904994.
18. Kim DN, Bathe KJ. A 4-node 3D-shell element to model shell surface tractions and incompressible behavior. *Comput Struct*. 2008;86(21-22):2027-2041. doi:10.1016/j.compstruc.2008.04.019.
19. Sussman T, Bathe KJ. 3D-shell elements for structures in large strains. *Comput Struct*. 2013;122:2-12. doi:10.1016/j.compstruc.2012.12.018.

20. Bathe K-J. Insights and Advances in the Analysis of Structures. In: *Proceedings Fifth International Conference on Structural Engineering, Mechanics and Computation -- SEMC 2013*. University of Cape Town: Taylor & Francis; 2013.
21. Hughes TJR, Carnoy E. Nonlinear finite element shell formulation accounting for large membrane strains. *Comput Methods Appl Mech Eng*. 1983;39(1):69-82. doi:10.1016/0045-7825(83)90074-9.
22. Hughes TJR, Liu WK. Nonlinear finite element analysis of shells: Part I. three-dimensional shells. *Comput Methods Appl Mech Eng*. 1981;26(3):331-362. doi:10.1016/0045-7825(81)90121-3.
23. Hughes TJR, Liu WK, Levit I. Nonlinear Dynamic Finite Element Analysis of Shells. *Nonlinear Finite Elem Anal Struct Mech*. 1981:151-168. doi:10.1007/978-3-642-81589-8_9.
24. *Ansys Theoretical*. Release 10. SAS IP; 2008.
25. Bathe K-J. *Finite Element Procedures*. New Jersey: Pearson Education; 1996.
26. Belytschko T, Schwer L. Large displacement, transient analysis of space frames. *Int J Numer Meth Eng*. 1977;11(August 1975):65-84.
27. Bonet J, Wood RD. *NONLINEAR CONTINUUM MECHANICS FOR FINITE ELEMENT ANALYSIS*. Cambridge university press; 2008.
28. Miller K, Joldes G, Lance D, Wittek A. Total Lagrangian explicit dynamics finite element algorithm for computing soft tissue deformation. *Commun Numer Methods Eng*. 2006;23(2):121-134. doi:10.1002/cnm.887.
29. Shantaram D, Owen DRJ, Zienkiewicz OC. Dynamic transient behaviour of two- and three-dimensional structures including plasticity, large deformation effects and fluid interaction. *Earthq Eng Struct Dyn*. 1976;4(6):561-578. doi:10.1002/eqe.4290040605.
30. Taylor RL, Govindjee S. *FEAP --A Finite Element Analysis Program Version 8.4 Parallel User Manual*. Berkeley: University of California; 2013.
31. Hughes TJR. *The Finite Element Method Linear Static and Dynamic Finite Element Analysis*. New York: Dover Publications; 2000.
32. Shabana AA. *Computational Continuum Mechanics*. Cambridge: Cambridge University Press; 2011. doi:10.1017/CBO9781139059992.
33. Liu S, Lin T. Elastkâ€□plastic dynamic analysis of structures using known elastic solutions. *Earthq Eng Struct …*. 1979;7(August 1978):147-159. http://onlinelibrary.wiley.com/doi/10.1002/eqe.4290070204/abstract.
34. Sathyamoorthy M. *No TitleNonlinear Analysis of Structures*. New York: CRC Press; 1997.
35. Macneal RH, Harder RL. A proposed standard set of problems to test finite element accuracy. *Finite Elem Anal Des*. 1985;1(1):3-20. doi:10.1016/0168-874X(85)90003-4.
36. Zienkiewicz OC. *The Finite Element Method*. third ed. London: McGraw-Hill; 1977.
37. *Comsol Multiphysics* ®. Version 4. Comsol; 2012.
38. MacNeal RH. *The NASTRAN Theoretical Manual.* Washington ;Athens Ga.: Scientific and Technical Information Office National Aeronautics and Space Administration ;;For sale by Computer Software Management and Information Center University of Georgia; 1972.
39. Pakravan A, Krysl P. Mean-strain 10-node tetrahedron with energy-sampling stabilization. *Int J Numer Methods Eng*. 2016. doi:10.1002/nme.5335.